\documentclass{article}

\usepackage{arxiv}

\usepackage[utf8]{inputenc} 
\usepackage[T1]{fontenc}    
\usepackage{hyperref}       
\usepackage{url}            
\usepackage{booktabs}       
\usepackage{amsfonts}       
\usepackage{nicefrac}       
\usepackage{microtype}      
\usepackage{lipsum}		
\usepackage{graphicx}
\usepackage{doi}

\usepackage{stmaryrd}
\usepackage{latexsym}
\usepackage{curves}
\usepackage{xcolor}
\usepackage{colortbl}
\usepackage{makecell}

\usepackage{mathrsfs}
\usepackage{amscd,amsfonts,amsmath,amssymb,amsthm,enumerate,epsfig,color,xcolor}
\usepackage{algorithm}
\usepackage{algpseudocode}

\begin{document}

\title{Partitioned neural network approximation for partial differential equations enhanced with Lagrange multipliers and localized loss functions}

\date{}

\author{Deok-Kyu Jang \\
		Department of Applied Mathematics\\ 
		Kyung Hee University\\ 
		Yongin, Republic of Korea \\
	\texttt{dkjang@khu.ac.kr} \\
	\And
	Kyungsoo Kim \\
		Department of Applied Mathematics\\ 
		Kyung Hee University\\ 
		Yongin, Republic of Korea \\
	\texttt{kyungsoo@khu.ac.kr} \\
	\And
	Hyea Hyun Kim$^*$ \\
		Department of Applied Mathematics\\ 
		Kyung Hee University\\ 
		Yongin, Republic of Korea \\
	\texttt{hhkim@khu.ac.kr}
}

\renewcommand{\shorttitle}{PNNA with Lagrange multipliers and localized loss functions}

\hypersetup{
pdftitle={Partitioned neural network approximation for partial differential equations enhanced with Lagrange multipliers and localized loss functions},
pdfsubject={math.NA, cs.LG, physics.comp-ph},
pdfauthor={Deok-Kyu Jang, Kyungsoo Kim, Hyea Hyun Kim},
pdfkeywords={partitioned neural networks, non-overlapping domain decomposition, Lagrange multipliers, localized loss function, iterative algorithm, partial differential equations},
}

\maketitle

\begin{abstract}
Partitioned neural network functions are used to approximate the solution of partial differential equations. The problem domain is partitioned into non-overlapping subdomains and the partitioned neural network functions are defined on the given non-overlapping subdomains. Each neural network function then approximates the solution in each subdomain. To obtain the convergent neural network solution, certain continuity conditions on the partitioned neural network functions across the subdomain interface need to be included in the loss function, that is used to train the parameters in the neural network functions. In our work, by introducing suitable interface values, the loss function is reformulated into a sum of localized loss functions and each localized loss function is used to train the corresponding local neural network parameters. In addition, to accelerate the neural network solution convergence, the localized loss function is enriched with an augmented Lagrangian term, where the interface condition and the boundary condition are enforced as constraints on the local solutions by using Lagrange multipliers. The local neural network parameters and Lagrange multipliers are then found by optimizing the localized loss function. To take the advantage of the localized loss function for the parallel computation,
an iterative algorithm is also proposed. For the proposed algorithms, their training performance and convergence are numerically studied for various test examples.
\end{abstract}

\section{Introduction}

Recently, there have been developed many successful approaches to solve partial differential equations (PDEs) using neural network approximation~\cite{sirignano2018,yu2018deep,raissi2019}.
The advantage of these new approaches is that they can be used for PDEs without any concern on the discretization method suitable for the given problem and the shape or dimension of the problem domain.
On the other hand, the performance and solution accuracy of these methods highly depend on hyperparameter settings, i.e., the neural network architecture, the choice of data sets for training the neural network parameters, the design of loss functions, and the optimization methods.
In general, even with nicely tuned hyperparameters they require long parameter training time
and  such new approximation approaches become thus impractical
compared to classical approximation methods, i.e., finite difference methods,
finite element methods, spectral element methods, and discontinuous Galerkin method,
when used for application problems.

To deal with such limitations of neural network approximation,
there have been some successful previous studies that utilize partitioned neural network functions
as a solution surrogate.
In the partitioned neural network approximation,
the problem domain is partitioned into many smaller subdomains and the partitioned neural network functions
are formed by local neural network functions, where each local neural network function approximates
the solution in each subdomain.
Such a partitioned structure allows more flexible network design and hyperparameter settings, and
it can also utilize parallel computation schemes to speed up the computation time
in parameter training process.
For problems with interface jumps inside the problem domain, with heterogeneous coefficients, or
with multi-physics model equations,
it is beneficial to use partitioned neural network functions as a solution surrogate
so as to handle such heterogeneities easily.
In \cite{cPINN,xPINN}, such an idea was first proposed for non-overlapping subdomain partitions of
the problem domain, where each partitioned neural network is trained for the localized model problem
in each subdomain with suitable continuity conditions across the subdomain interface and
the methods showed promising results for more challenging application models, that could not be well
resolved with a single larger neural network function.
In \cite{FBPINN}, a similar approach was proposed for overlapping subdomain partitions by
forming a neural network solution as a sum of localized neural network functions,
where the neural network functions are localized in each subdomain by using window functions
and their parameters are trained for the global problem without the need for the additional continuity condition, in contrast to the non-overlapping subdomain case~\cite{cPINN,xPINN}.
We note that the localized nature of parameters in the partitioned neural network functions
has not been fully utilized when training parameters in the previous studies~\cite{cPINN,xPINN,FBPINN}.
Every training epoch, the localized parameters are updated with the corresponding gradient value of the loss function. The evaluation of the gradient value needs data communication, i.e., the exchange of parameter values, between the neighboring subdomains.
Since the number of training epochs is often more than several tens of thousands in many application problems, such data communication becomes a bottleneck for parallel computation.

To alleviate such data communication issue in the partitioned neural network approximation,
iteration algorithms in domain decomposition methods~\cite{TW-Book} can be adopted for training the local neural network functions.
In the iteration algorithms, the data communication is needed every outer iteration
and inside each outer iteration the local neural network is trained for the local problem with
a provided interface condition.
In \cite{DD26-kim-yang,DD27-kim-yang}, iteration methods based on additive Schwarz algorithms have been developed
for the neural network approximation, that is formed for the overlapping subdomain partitions, and the methods
showed efficient training results under the parallel computing environments.
Rigorous convergence analysis and a faster two-level iteration algorithm have been
also studied in a more recent work~\cite{arXiv-yang-kim-2023}.
We also note that there have been some previous similar studies~\cite{li2019,li2020-pro}
to address the parallel computation issue in the neural network approximation by utilizing domain decomposition
algorithms.

The purpose of this work is to extend the third author's previous study~\cite{DD27-kim-yang} to neural network approximation with non-overlapping subdomain partitions.
Differently to the overlapping subdomain partition case,
additional continuity conditions on the local neural network solutions
across the subdomain interface are needed and they are included to
the loss function for training the parameters in the local neural network solutions.
In our loss function development, we introduce solution and flux values on the subdomain interface
and reformulate the standard loss function used in \cite{cPINN} as a sum of localized loss functions.
Our resulting localized loss function in each subdomain consists of several loss terms, i.e.,
the PDE residual loss, the boundary value loss, and the boundary normal flux loss.
For such a case, the imbalance among different loss terms can often cause a failure in the parameter training
and carefully chosen weight factors are needed for the success of the parameter training.
However, such hyperparameter tunings are difficult in practice.
To deal with this issue, we include an augmented Lagrangian relaxation term to the localized loss function.
The augmented Lagrangian term is formed by enforcing the boundary condition for the local solution
as constraints with Lagrange multipliers.
We note that in \cite{Son2023} a method of an augmented Lagrangian relaxation for the boundary condition
was proposed to enhance the training efficiency and accuracy of neural network solutions for differential equations.
A similar approach was also developed for variational problems with essential boundary conditions in \cite{Huang2021}.
However, such approaches have not been applied to partitioned neural network approximation.
In \cite{Scalable-PINN2022}, Lagrange multipliers are introduced
as trainable self-adaptation weights to the loss function, where Lagrange multipliers
are optimized to reduce the corresponding loss term of data points.
We note that in \cite{Scalable-PINN2022} partitioned neural network approximation is considered but their use of Lagrange multipliers is a non-standard way and is different from those in \cite{Son2023,Huang2021}.
As our final goal, taking the advantage of our localized loss functions,
we propose an iteration algorithm so as to enhance the training efficiency further
under the parallel computing environment.
Concrete convergence analysis and further acceleration schemes
for our proposed iteration algorithm
will be studied in our forthcoming research.

This paper is organized as follows. In Section~\ref{sec:methods}, we develop partitioned neural network approximation methods for solving Poisson model problems by introducing localized loss functions and augmented Lagrangian terms.
Three algorithms are proposed and their convergence and training performance are numerically studied in Section~\ref{sec:numerics}.
The proposed algorithms are further applied to more challenging examples of a discontinuous coefficient elliptic problem, a Poisson interface problem, and a Stokes flow with an immersed interface.
In Section~\ref{sec:conclude}, we summarize our proposed methods and discuss some future works for improvements.

\section{Partitioned neural network approximation with localized loss functions, augmented Lagrangian relaxation, and an iteration method}\label{sec:methods}
In this section, we propose three algorithms (Algorithms~\ref{algo1}-\ref{algo3})
for the partitioned neural network approximation
to PDEs.
The neural network approximation is formed by localized neural network functions
defined on a non-overlapping subdomain partition of the problem domain.
Each local neural network function approximates the solution on the corresponding subdomain in the partition.

The parameter in each local neural network function is trained for the loss function,
that is formed so as to satisfy the localized PDEs, boundary conditions, and appropriate interface conditions
across the subdomain boundary.
The appropriate interface conditions are crucial in obtaining
convergent neural network approximation to the problem solution.
In Algorithm~\ref{algo1},
additional unknowns for the solution and the flux values on the subdomain interface
will be introduced and using them the loss function in \cite{cPINN,xPINN}
is reformulated as a sum of localized loss functions.
The localized loss functions are easily adoptable for the parallel computation than
those used in \cite{cPINN,xPINN}.
In Algorithm~\ref{algo2}, the localized loss functions
will be further enriched with augmented Lagrangian relaxation terms
so as to enhance the approximate solution accuracy and the parameter training performance.
Finally, an iterative algorithm, i.e., Algorithm~\ref{algo3}, will be proposed for training
parameters in the partitioned neural network functions, that will require much lesser data communication
between neighboring subdomains than those in Algorithms~\ref{algo1} and \ref{algo2}.

Our method development will be presented for the Poisson model problem in a two-dimensional domain $\Omega$.
We note that our methods can be directly extended to more general second order PDEs
by introducing suitable continuity conditions on the subdomain interface.
In our numerical experiments later, we also extend the proposed methods to a Poisson problem with discontinuous coefficients, an interface problem, and the Stokes problem with an immersed interface
to show the performance of the proposed methods.

\subsection{Algorithm 1: Partitioned neural network approximation with localized loss functions}
For the method description, we consider the following Poisson problem in a two-dimensional domain $\Omega$,
\begin{equation}\label{model:Poisson}
\begin{split}
-\triangle u({\bf x})&=f({\bf x}),\quad \forall {\bf x} \in \Omega,\\
u({\bf x})&=g({\bf x}),\quad \forall {\bf x } \in \partial \Omega,
\end{split}
\end{equation}
where we assume that the solution of the above model problem uniquely exists.
By introducing a non-overlapping subdomain partition $\{ \Omega_i \}$ of the domain $\Omega$, the following local
Poisson problem in each subdomain $\Omega_i$ can be considered,
\begin{equation}\label{model:local:Poisson}
\begin{split}
-\triangle u_i({\bf x})&=f({\bf x}),\quad \forall {\bf x} \in \Omega_i,\\
u_i({\bf x})&=g({\bf x}),\quad \forall {\bf x} \in \partial \Omega_i \bigcap \partial \Omega.
\end{split}
\end{equation}
In order to obtain the local problem solution $u_i$ as the restriction of the solution $u$ to each subdomain $\Omega_i$,
the following additional continuity conditions for $u_i$ on the subdomain interface are needed,
\begin{align}
u_i({\bf x})&=u_j({\bf x}),\quad \forall {\bf x} \in \Gamma_{ij},\label{u:conti}\\
\frac{\partial u_i}{\partial n_i}({\bf x})&=-\frac{\partial u_j}{\partial n_j}({\bf x}),\quad \forall {\bf x} \in \Gamma_{ij},\label{un:conti}
\end{align}
where $\Gamma_{ij}$ denotes the interface of the two subdomains $\Omega_i$ and $\Omega_j$,
and $\frac{\partial u_i}{\partial n_i}$ and $\frac{\partial u_j}{\partial n_j}$ denote the derivative in the outward normal direction on the subdomain interface.

In order to include more general interface problems, we can rewrite the above local problems into a more general form,
\begin{equation}\label{eq_poi_dd}
\begin{split}
-\triangle u_i &=f \quad  \text{in } \Omega_i, \\
u_i &= g \quad \text{on } \Gamma_{i0} = \partial \Omega_i \bigcap \partial \Omega, \\
\left\llbracket  u\right\rrbracket_{\Gamma_{ij}} &= p_{ij}  \quad \text{on } \Gamma_{ij}, \ \forall j \in s(i), \\
\left\llbracket  \frac{\partial u}{\partial n} \right\rrbracket _{\Gamma_{ij}} &= q_{ij}  \quad \text{on } \Gamma_{ij}, \ \forall j \in s(i),
\end{split}
\end{equation}
where the notations $\left\llbracket  u\right\rrbracket_{\Gamma_{ij}}$ and
$\left\llbracket  \frac{\partial u}{\partial n } \right\rrbracket_{\Gamma_{ij}}$ are defined as
$$\left\llbracket  u\right\rrbracket_{\Gamma_{ij}}=(u_i {\bf n}_i + u_j {\bf n}_j) \cdot {\bf n},\quad
\left\llbracket  \frac{\partial u}{\partial {n}} \right\rrbracket_{\Gamma_{ij}}
=\left(\frac{\partial u_i}{\partial n_i} {\bf n}_i +\frac{\partial u_j}{\partial n_i} {\bf n}_j \right)\cdot {\bf n} ,$$
with ${\bf n}$ being a fixed unit normal vector to the interface $\Gamma_{ij}$.
We note that there are ${\bf n}_i$ and ${\bf n}_j$ unit normal vectors to the interface $\Gamma_{ij}$ and
each corresponds to the outward normal vector to the corresponding subdomain.
In addition, $p_{ij}$ and $q_{ij}$ are the given jump conditions for the local solutions on the interface, and $s(i)$ denotes the set of subdomain indices $j$ such that $\Omega_j$ is the neighboring subdomain sharing the interface $\Gamma_{ij}$.
We note that for the Poisson problem~\eqref{model:Poisson}
the interface jump conditions vanish, i.e., $p_{ij}=0$ and $q_{ij}=0$.

After we partition the model problem~\eqref{model:Poisson} into each subdomain problem as in the above~\eqref{eq_poi_dd}, we now proceed to find an approximate solution by using partitioned neural network functions.
The advantages in using the partitioned problems are first, the partitioned neural networks can be used as an approximate solution to allow more flexible hyperparameter settings, and second, training local neural network parameters in the partitioned neural networks can be done more efficiently than training the parameters in the single large neural network function defined on the whole problem domain.

We now introduce a local neural network function $U_i({\bf x}; \theta_i)$, that approximates the local problem solution $u_i$ in \eqref{eq_poi_dd}, and choose the training data sets $X_{\Gamma_{ij}}$, $X_{\Gamma_{i0}}$,
and $X_{\Omega_i}$ from each corresponding interfaces and each subdomain, i.e., $\Gamma_{ij}$, $\Gamma_{i0}$,
and $\Omega_i$.
We construct required loss terms to train $U_i$, following the approaches in \cite{raissi2019},
\begin{equation}\label{loss1}
\mathcal{L}_{f,i} = \frac{1}{\vert X_{\Omega_i} \vert} \sum_{{\bf x} \in X_{\Omega_i}} \vert \triangle U_i ({\bf x}; \theta_i) + f({\bf x}) \vert^2,  \qquad
\mathcal{L}_{g,i} = \frac{1}{\vert X_{\Gamma_{i0}} \vert} \sum_{{\bf x} \in X_{\Gamma_{i0}}} \vert U_i ({\bf x}; \theta_i) - g({\bf x}) \vert^2.
\end{equation}
The loss terms listed in \eqref{loss1} are formed so as to satisfy the given PDE and boundary conditions. To satisfy
the condition on the solution and its normal derivative across the subdomain interface, see \eqref{eq_poi_dd}, we need to introduce the following additional loss terms as in \cite{cPINN,xPINN}:
\begin{equation}\label{loss2_1}
\begin{split}
&\sum_{j \in s(i)} \frac{1}{\vert X_{\Gamma_{ij}} \vert} \sum_{{\bf x} \in X_{\Gamma_{ij}} } \vert (U_i ( {\bf x}; \theta_i){\bf n}_i + U_j ( {\bf x}; \theta_j) {\bf n}_j) \cdot {\bf n} - p_{ij}( {\bf x} ) \vert^2,  \\
&\sum_{j \in s(i)} \frac{1}{\vert X_{\Gamma_{ij}} \vert} \sum_{{\bf x} \in X_{\Gamma_{ij}} }  \left\vert  \left( \frac{\partial U_i ({\bf x}; \theta_i)}{\partial n_i} {\bf n}_i + \frac{\partial U_j ({\bf x}; \theta_i)}{\partial n_i} {\bf n}_j \right) \cdot {\bf n} - q_{ij}({\bf x}) \right\vert^2.
\end{split}
\end{equation}
We note that the sum of the above mentioned loss functions in \eqref{loss1}-\eqref{loss2_1} is used to train the parameters
$\theta_i$ for the local neural network solutions $U_i({\bf x};\theta_i)$ in the previous study~\cite{cPINN}.
However, the resulting loss function is then coupled on the subdomain interface.

To achieve a more stable optimization performance and to localize the loss construction for each $U_i$,
the loss terms in \eqref{loss2_1} are reformulated as follows.
To make the presentation below simpler, we set ${\bf n}={\bf n}_i$ for $i>j$ on the interface $\Gamma_{ij}$.
We first define $\tilde{U}_{ij}({\bf x})$, $\tilde{U}_{n,ij}({\bf x})$,
$\tilde{U}_{ji}({\bf x})$, and $\tilde{U}_{n,ji}({\bf x})$ for $i>j$ as
\begin{equation}\label{u_tildeij}
\begin{split}
\tilde{U}_{ij} ({\bf x}) &= \frac{1}{2}( U_i ({\bf x}; \theta_i) + U_j ({\bf x}; \theta_j))+\frac{1}{2} p_{ij}({\bf x}), \\
\tilde{U}_{n,ij}({\bf x})&= \frac{1}{2} \left( \frac{\partial U_i ({\bf x}; \theta_i)}{\partial n} + \frac{\partial U_j ({\bf x}; \theta_j)}{\partial n}  \right) + \frac{1}{2} q_{ij}({\bf x})
\end{split}
\end{equation}
and
\begin{equation}\label{u_tildeji}
\begin{split}
\tilde{U}_{ji} ({\bf x}) &= \frac{1}{2} ( U_i ({\bf x}; \theta_i) + U_j ({\bf x}; \theta_j) )-\frac{1}{2} p_{ij}({\bf x}), \\
\tilde{U}_{n,ji}({\bf x}) &= \frac{1}{2} \left(\frac{\partial U_i ({\bf x}; \theta_i)}{\partial n} + \frac{\partial U_j ({\bf x}; \theta_j)}{\partial n}   \right) - \frac{1}{2} q_{ij}({\bf x}).
\end{split}
\end{equation}
We note that when $p_{ij}=0$ the definitions for $\tilde{U}_{ij}({\bf x})$ and $\tilde{U}_{ji}({\bf x})$
are identical. The same holds for $\tilde{U}_{n,ij}$ and $\tilde{U}_{n,ji}$ when $q_{ij}=0$.
By using \eqref{u_tildeij} and \eqref{u_tildeji}, we rewrite the jump condition in \eqref{loss2_1} into
\begin{equation}\label{Uij:split}
\begin{split}
(U_i ({\bf x}; \theta_i){\bf n}_i +  U_j ({\bf x}; \theta_j){\bf n}_j) \cdot {\bf n} - p_{ij}({\bf x}) &=
U_i({\bf x};\theta_i)-U_j({\bf x};\theta_j)-p_{ij}({\bf x})\\
&=(U_i({\bf x};\theta_i)-\tilde{U}_{ij})-(U_j({\bf x};\theta_j)-\tilde{U}_{ji})
\end{split}
\end{equation}
and
\begin{equation}\label{Unij:split}
\begin{split}
\left( \frac{\partial U_i ({\bf x}; \theta_i)}{\partial n_i} {\bf n}_i + \frac{\partial U_j ({\bf x}; \theta_i)}{\partial n_i} {\bf n}_j \right) \cdot {\bf n} - q_{ij}({\bf x})&=
\left( \frac{\partial U_i({\bf x};\theta_i)}{\partial n}-\frac{\partial U_j({\bf x};\theta_j)}{\partial n} \right)-q_{ij}({\bf x})\\
&=\left( \frac{\partial U_i ({\bf x};\theta_i)}{\partial n}-\tilde{U}_{n,ij}({\bf x})\right)
-\left(\frac{\partial U_j({\bf x};\theta_j)}{\partial n}-\tilde{U}_{n,ji}({\bf x})\right)
\end{split}
\end{equation}
to define the following localized loss terms for the interface conditions,
\begin{equation}\label{loss2}
\begin{split}
\mathcal{F}_{u,i} &= \sum_{j \in s(i)} \frac{1}{\vert X_{\Gamma_{ij}} \vert} \sum_{{\bf x} \in X_{\Gamma_{ij}} } \vert U_i ({\bf x}; \theta_i) - \tilde{U}_{ij} ({\bf x}) \vert^2,  \\
\mathcal{F}_{n,i} &= \sum_{j \in s(i)} \frac{1}{\vert X_{\Gamma_{ij}} \vert} \sum_{{\bf x} \in X_{\Gamma_{ij}} } \left\vert  \frac{\partial U_i ({\bf x}; \theta_i)}{\partial n} - \tilde{U}_{n,ij} ({\bf x}) \right\vert^2.
\end{split}
\end{equation}
By using the loss terms in \eqref{loss2},
we can then form the following minimization problem for $U_i ({\bf x}, \theta_i)$;
find the parameter $\theta_i$ such that
\begin{equation}\label{loss_algo1}
\min_{\theta_i} \mathcal{J}_i (\theta_i),
\end{equation}
where
\begin{equation}\label{loss_algo1_terms}
\mathcal{J}_i (\theta_i) :=  \mathcal{L}_{f,i} + \mathcal{L}_{g,i} + \mathcal{F}_{u,i} + \mathcal{F}_{n,i}.
\end{equation}
Since the terms $\tilde{U}_{ij} ({\bf x})$ and $\tilde{U}_{n,ij}({\bf x})$ contain
both parameters $\theta_i$ and $\theta_j$,
we propose the following Algorithm~\ref{algo1} for training the parameters $\theta_i$, where the values $\tilde{U}_{ij} ({\bf x})$ and $\tilde{U}_{n,ij}({\bf x})$ are computed
with $\theta_i^{(\ell)}$ and $\theta_j^{(\ell)}$, i.e., the trained parameters from the previous epoch, and the parameter $\theta_i^{(\ell+1)}$ is then updated by using the localized loss terms with the provided interface conditions $\tilde{U}_{ij}({\bf x})$ and $\tilde{U}_{n,ij}({\bf x})$.
The detailed description of the parameter training procedure
with the localized loss functions is listed in Algorithm~\ref{algo1}.

\begin{algorithm}[ht!]
\caption{PNNA-LL (Partitioned Neural Network Approximation with Localized Loss functions)}\label{algo1}
\begin{algorithmic}
\Require $N$=Maximum epochs
\State Set initial $\theta_i$ for all $i$ and set $\theta_i^* := \theta_i$.
\State For all $i$, calculate initial $\tilde{U}_{ij}$ and $\tilde{U}_{n,ij}$  for all $s(i)$.
\State $\ell=0$
\While{iteration : $\ell < N$}
	\State Update $\theta_i$ by minimizing $\mathcal{J}_i :=  \mathcal{L}_{f,i} + \mathcal{L}_{g,i} + \mathcal{F}_{u,i} + \mathcal{F}_{n,i}$ for all $i$.
	\State  {\bf Communication} : For all $i$, update $\tilde{U}_{ij}$ and $\tilde{U}_{n,ij}$ for all $s(i)$
	\If{$ \mathcal{J}_i (\theta_i) < \mathcal{J}_i (\theta_i^*)$ for each $i$}
			\State $\theta_i^*$ = $\theta_i$
	\EndIf
    \State $\ell=\ell+1$
\EndWhile
\end{algorithmic}
\end{algorithm}

In all the proposed algorithms of our work, as the final output we take the parameter $\theta_i^*$ that corresponds to the minimum loss value during the $N$ training epochs.
The parameter $\theta_i^*$ is updated whenever the minimum loss value
occurs during the parameter training.
Algorithm~\ref{algo1} is similar to those already studied in \cite{cPINN,xPINN}.
One advantage in Algorithm~\ref{algo1} is that the local neural network parameter
$\theta_i$ is trained for the localized loss function.
Based on the localized loss function, we can further improve the proposed algorithm to be more suitable
for parallel computing environments.
For that purpose, we will propose an iterative algorithm where each local problem is solved by training
the local neural network solution for the given iterate and the iterate is updated with
the trained local solutions to proceed the next iteration so as to obtain
a convergent neural network approximate solution.

\subsection{Algorithm 2: Localized loss function enriched with augmented Lagrangian terms}
Before we develop the iterative algorithm, we improve the localized loss function in \eqref{loss_algo1_terms}
by including augmented Lagrangian terms for the interface conditions, inspired by the work in \cite{Son2023},
where the use of Lagrange multipliers to the boundary conditions is shown to be effective for
improving the solution accuracy and training performance in the neural network approximation.
Their work was studied only for a single neural network case,
that approximates the solution in the whole problem domain.
In our Algorithm~\ref{algo2}, we extend their method
to the partitioned neural network approximation.

We recall the loss function in \eqref{loss_algo1_terms}. Among the loss terms, those related
to the boundary condition $U_i=g$ on $\Gamma_{i0}$ and the interface condition $U_i=\widetilde{U}_{ij}$ on $\Gamma_{ij}$ are more difficult to be trained than the other terms and they need specially tuned hyperparameter settings, such as large weight factors to penalize the corresponding loss terms.
As an alternative to the hyperparameter tunings,
the boundary condition in $\mathcal{L}_{g,i}$ and the interface condition in $\mathcal{F}_{u,i}$
can be enforced as constraints on the local solution $U_i({\bf x};\theta_i)$ by using Lagrange multipliers,
similarly as in \cite{Son2023}.
For that, we now introduce vector unknowns $\lambda_{ij}$ such that their lengths are identical to the number of points in the training data set $X_{\Gamma_{ij}}$.
For $\Gamma_{i0}$, we introduce the vector unknown $\lambda_{i0}$ similarly.
With the Lagrange multipliers $\lambda_{ij}$ and $\lambda_{i0}$, the constraints on the local solution $U_i({\bf x},\theta_i)$ are implemented by introducing the following terms,
\begin{equation}\label{loss3}
\begin{split}
\Lambda_{i0}(\theta_i,\lambda_{i0}) &=  \sum_{{\bf x} \in X_{\Gamma_{i0}} } \lambda_{i0}({\bf x}) ( U_i ({\bf x}; \theta_i) - g(x)),\\
\Lambda_{i}(\theta_i,(\lambda_{ij})_{j}) &=  \sum_{j \in s(i)} \sum_{{\bf x} \in X_{\Gamma_{ij}} }  \lambda_{ij}({\bf x}) ( U_i ({\bf x}; \theta_i) - \tilde{U}_{ij} ({\bf x}) ),
\end{split}
\end{equation}
where the notation $(\lambda_{ij})_j$ denotes the union of all the Lagrange multipliers $\lambda_{ij}({\bf x})$, i.e.,
$$(\lambda_{ij})_j:=\{ \lambda_{ij}({\bf x})\,;\, \forall {\bf x} \in X_{\Gamma_{ij}},\; \forall j \in s(i) \}.$$
By including the additional terms into the loss function,
we can now form the following min-max problem for $U_i ({\bf x}; \theta_i)$, $\lambda_{i0}$, and $(\lambda_{ij})_{j}$; find the parameters $\theta_i$, $\lambda_{i0}$, and $(\lambda_{ij})_{j}$  such that
$$ \max_{(\lambda_{ij})_{j}, \lambda_{i0}} \min_{\theta_i} \mathcal{J}_{i,\Lambda}(\theta_i,\lambda_{i0},(\lambda_{ij})_j )$$
where the loss function is defied as
\begin{equation}\label{loss:min-max}
\mathcal{J}_{i,\Lambda}(\theta_i,\lambda_{i0},(\lambda_{ij})_j ):=\mathcal{J}_i (\theta_i)+\Lambda_{i0}(\theta_i,\lambda_{i0})+\Lambda_i(\theta_i,(\lambda_{ij})_{j}).
\end{equation}
We note that definitions for $\mathcal{J}_i(\theta_i)$, $\Lambda_{i0}(\theta_i,\lambda_{i0})$, and $\Lambda_i(\theta_i,(\lambda_{ij})_{j})$
are given in \eqref{loss_algo1_terms} and \eqref{loss3}.
In the above min-max problem, Lagrange multipliers $\lambda_{ij}$ and $\lambda_{i0}$ are included to enforce the constraints in the partitioned problem and they are updated to maximize the loss function in \eqref{loss:min-max}.
By using the gradient based method, every training epoch the Lagrange multipliers are then updated
as
\begin{equation}\label{LM:update}
\begin{split}
\lambda_{i0}^{(\ell)} ({\bf x}) &= \lambda_{i0}^{(\ell-1)} ({\bf x}) + \alpha_0 \nabla_{\lambda_{i0}} \mathcal{J}_{i,\Lambda}=\lambda_{i0}^{(\ell-1)} ({\bf x}) + \alpha_0 ( U_i ({\bf x}; \theta_i^{(\ell)}) - g ({\bf x}) ),\\
\lambda_{ij}^{(\ell)} ({\bf x}) &= \lambda_{ij}^{(\ell-1)} ({\bf x}) + \alpha_{\lambda} \nabla_{\lambda_{ij}} \mathcal{J}_{i,\Lambda} =\lambda_{ij}^{(\ell-1)} ({\bf x}) + \alpha_{\lambda}( U_i ({\bf x}; \theta_i^{(\ell)}) - \tilde{U}_{ij} ({\bf x}) ),
\end{split}
\end{equation}
where $\ell$ denotes the training epochs, and $\alpha_{\lambda}$ and $\alpha_0$ are the learning rates.
The detailed description of the parameter training procedure for $\theta_i$, $\lambda_{i0}$, and $\lambda_{ij}$
is provided in Algorithm~\ref{algo2}.

\begin{algorithm}[ht!]
\caption{PNNA-LL-AL (PNNA-LL with Augmented Lagrangian terms)}\label{algo2}
\begin{algorithmic}
\Require $N$=Maximum epochs, learning rates $\alpha_\lambda$ and $\alpha_0$
\State Set initial $\theta_i$ for all $i$ and set $\theta_i^* := \theta_i$.
\State For all $i$, set initial $\lambda_{i0}$ and $\lambda_{ij}$ for all $j \in s(i)$.
\State For all $i$, calculate initial $\tilde{U}_{ij}$ and $\tilde{U}_{n,ij}$  for all $j \in s(i)$.
\State $\ell=0$
\While{iteration : $\ell<N$}
	\State Update $\theta_i$ by minimizing $\mathcal{J}_{i,\Lambda}  =  \mathcal{J}_i + \Lambda_{i0} + \Lambda_{i}$ for all $i$. $\left( {\mathcal{J}}_i :=  \mathcal{L}_{f,i} + \mathcal{L}_{g,i} + \mathcal{F}_{u,i} + \mathcal{F}_{n,i} \right)$
	\State {\bf Communication} : For all $i$, update $\tilde{U}_{ij}$ and $\tilde{U}_{n,ij}$ for all $s(i)$
	\State Update $\lambda_{i0}$ and $\lambda_{ij}$ by gradient ascent.
	\If{$ {\mathcal{J}}_i (\theta_i) < {\mathcal{J}}_i (\theta_i^*)$ for each $i$,}
			\State $\theta_i^*$ = $\theta_i$
	\EndIf
    \State $\ell=\ell+1$
\EndWhile
\end{algorithmic}
\end{algorithm}

We would like to stress that the inclusion of the augmented Lagrangian term
helps to increase the solution accuracy and training performance compared to the standard loss function case, as seen in the numerical results listed
in Table~\ref{sec31_tb1} for the single neural network approximation.
In our development of the iteration method in the next subsection, we will thus consider Algorithm~\ref{algo2}.
The proposed iteration method can be applied directly to Algorithm~\ref{algo1} as well.

\subsection{Algorithm 3: Iteration method for partitioned neural network approximation}

Algorithm~\ref{algo2} should update $\theta_i$ and perform communication between adjacent subdomains using the updated neural network solutions in order to calculate $\tilde{U}_{ij}({\bf x})$ and $\tilde{U}_{n,ij}({\bf x})$.
Except that, the gradient update for $\theta_i$, $\lambda_{i0}$, and $\lambda_{ij}$ can be done independently
for the localized loss function.
To take advantage of such a localized feature of the parameter training procedure,
we now propose an iteration method for Algorithm~\ref{algo2}, where the interface values $\tilde{U}_{ij}({\bf x})$ and $\tilde{U}_{n,ij}({\bf x})$ are updated after training the local neural network solutions
for the given interface condition.

Our iteration algorithm is based on the iteration methods developed in classical domain decomposition methods~\cite{TW-Book}.
In the classical domain decomposition iterative algorithms,
when the domain is partitioned into non-overlapping subdomains, the interface solution values
or the flux values are updated after solving local problems independently for the interface values provided by the previous iterate.
The next iteration is proceeded with the updated interface values repeatedly until the solution convergence.
We also note that many successful domain decomposition algorithms with similar update
procedures on the interface solution values have been developed in~\cite{Substruct,FETI,BDD} for non-overlapping
subdomain partitions.

The description of our iteration method is listed in Algorithm~\ref{algo3}, where for the given initial
iterates $\tilde{U}_{ij}^{(0)}$ and $\tilde{U}_{n,ij}^{(0)}$ the local problems
are solved with the localized loss function in which $\tilde{U}_{ij}$ and $\tilde{U}_{n,ij}$
are set to $\tilde{U}_{ij}^{(0)}$ and $\tilde{U}_{n,ij}^{(0)}$.
When solving the local problems, we train the local parameters $\theta_i$, $\lambda_{i0}$, and $\lambda_{ij}$
to optimize the localized loss function and the parameters are updated by using the gradient method
up to $N_l$ number of training epochs. After that, the trained parameters are set as $\theta_i^{(1)}$.
They are used to update the iterates $\tilde{U}_{ij}^{(1)}$ and $\tilde{U}_{n,ij}^{(1)}$,
$$\tilde{U}_{ij}^{(1)}({\bf x})=\frac{1}{2}\left(U_i({\bf x};\theta_i^{(1)})+U_j({\bf x};\theta_j^{(1)})\right)+\frac{1}{2}p_{ij}({\bf x}), \; i>j,$$
$$\tilde{U}_{ij}^{(1)}({\bf x})=\frac{1}{2}\left(U_i({\bf x};\theta_i^{(1)})+U_j({\bf x};\theta_j^{(1)})\right)-\frac{1}{2}p_{ij}({\bf x}), \;
i<j,$$
and
$$\tilde{U}_{n,ij}^{(1)}({\bf x})=\frac{1}{2}\left(\frac{\partial U_i({\bf x};\theta_i^{(1)})}{\partial n} + \frac{\partial U_j({\bf x};\theta_j^{(1)})}{\partial n}\right)+\frac{1}{2} q_{ij}({\bf x}),\; i>j,$$
$$\tilde{U}_{n,ij}^{(1)}({\bf x})=\frac{1}{2}\left(\frac{\partial U_i({\bf x};\theta_i^{(1)})}{\partial n} + \frac{\partial U_j({\bf x};\theta_j^{(1)})}{\partial n}\right)-\frac{1}{2} q_{ij}({\bf x}),\; i<j.$$
In the next iteration, the updated values are used for the interface conditions, $\tilde{U}_{ij}$ and
$\tilde{U}_{n,ij}$ in the localized loss function when training the local neural network solutions.
Such an outer iteration is then performed up to $N_o$ number of iterations, with $N_o$
much smaller than $N$, the maximum epochs in Algorithms~\ref{algo1} and \ref{algo2}.
We note that in Algorithm~\ref{algo3}, the communication between the neighboring local neural network solutions is needed once before starting the next outer iteration after obtaining the local neural network solutions
for the given current interface values, $\tilde{U}_{ij}$ and $\tilde{U}_{n,ij}$.
The communication cost in Algorithm~\ref{algo3} is thus greatly reduced compared to the previous Algorithms~\ref{algo1} and \ref{algo2}.

\begin{algorithm}[ht!]
\caption{(Iteration Algorithm for PNNA-LL-AL)}\label{algo3}
\begin{algorithmic}
\Require $N_o$=Maximum outer iterations, $N_l$=Maximum epochs for local problem, learning rates $\alpha_\lambda$ and $\alpha_0$
\State Set initial $\theta_i$ for all $i$ and set $\theta_i^* := \theta_i$.
\State For all $i$, set initial $\lambda_{i0}$.
\State For all $i$, calculate initial $\tilde{U}_{ij}$ and $\tilde{U}_{n,ij}$  for all $j \in s(i)$.
\State $\ell_o=0$
\While{Outer iteration : $\ell_o < N_o$}
	\State For all $i$, set initial $\lambda_{ij}$ for all $j \in s(i)$.
	\State $\ell=0$
	\While {Local problem training epochs : $\ell < N_{l}$}
		\State Update $\theta_i$ by minimizing $\mathcal{J}_{i,\Lambda}  =  \mathcal{J}_i + \Lambda_{i0} + \Lambda_{i}$ for all $i$. $\left( {\mathcal{J}}_i :=  \mathcal{L}_{f,i} + \mathcal{L}_{g,i} + \mathcal{F}_{u,i} + \mathcal{F}_{n,i} \right)$
		\State Update $\lambda_{i0}$ and $\lambda_{ij}$ by gradient ascent.
		\If{$ {\mathcal{J}}_i (\theta_i) < {\mathcal{J}}_i (\theta_i^*)$ for each $i$,    }
				\State $\theta_i^*$ = $\theta_i$
		\EndIf
        \State $\ell=\ell+1$
	\EndWhile
	\State {\bf Communication} : For all $i$, update $\tilde{U}_{ij}$ and $\tilde{U}_{n,ij}$ for all $s(i)$
    \State $\ell_o=\ell_o+1$
\EndWhile
\end{algorithmic}
\end{algorithm}

In the following section, we will present the performance of the proposed Algorithms~\ref{algo1}-\ref{algo3}
for various test examples.
The theoretical convergence proof for the proposed Algorithms~\ref{algo1}-\ref{algo3} will be provided
in our forthcoming research.
Interestingly Algorithm~\ref{algo3} is closely related to the localized finite element tearing and interconnecting (FETI) methods, see~\cite{local-FETI,total-FETI}.
Based on this observation, we also anticipate that
the convergence of the proposed algorithms can be analyzed by using
well-developed analytical tools in the FETI methods, see~\cite{TW-Book,FETI}.


\section{Numerical results}\label{sec:numerics}
In our numerical experiments, we use the following settings for the neural networks and training their parameters.
Depending on the model problem and the subdomain partition, we choose the number of hidden layers, the number of nodes per hidden layer, and the training data sets differently.
Those specific details are provided later.
For all the experiments, we choose $\sin(x)$ as the activation function in the neural networks.
We note that the $\tanh(x)$ is a commonly used activation function in physics-informed neural networks (PINNs), however, our choice, the sine activation function was observed to be more effective and suitable for approximating derivative values, see \cite{sitzmann2020implicit}.
Network parameters are initialized by Xavier normal initialization (gain$=1.0$) for weights and by the initial biases
filled with 0.01.
To train the neural network parameters, we use the Adam optimizer~\cite{adam} with a fixed learning rate $lr=0.001$.
All samples are distributed by the Latin Hypercube sampling method.
We employ the following three algorithms introduced in the previous section;
\begin{itemize}
\item Algorithm \ref{algo1} (denoted $\mathcal{A}_1$) : The algorithm needs data communication at each training epoch and does not use the additional Lagrange multipliers for the constraints.
\item Algorithm \ref{algo2} (denoted $\mathcal{A}_2$) : The algorithm needs data communication at each training epoch and uses the Lagrange multipliers for the boundary and interface constraints. All Lagrange multipliers are initialized to the zero value and updated every training epoch by the gradient ascent method.
\item Algorithm \ref{algo3} (denoted $\mathcal{A}_3$) : The algorithm needs data communication at each outer iteration, i.e., every $N_l$ training epochs. At the first outer iteration, when training the local neural network parameters, the Lagrange multipliers are set to zero and updated every training epoch by the gradient ascent method. After then, in the next outer iteration, the Lagrange multipliers for the interface constraints are initialized to the gradient value of the loss function $\nabla_{\lambda_{ij}}\mathcal{J}_{i,\Lambda}$, i.e.,
    $$\lambda_{ij}^{(0)} ({\bf x}) =  U_i ({\bf x}; \theta_i) - \tilde{U}_{ij} ({\bf x}) ,$$
 where $\theta_i$ are the trained parameters from the previous outer iteration and the values $\tilde{U}_{ij}$ are updated by the data communication between the trained neighboring neural network solutions.
 They are then updated by the gradient ascent method every training epoch of local problems up to the maximum training epochs $N_l$.
\end{itemize}
To estimate the performance of our proposed methods, we define the relative $L^2$-error $\epsilon_u$ as
\begin{equation}\label{rel_L2_err}
\epsilon_u(X)= \sqrt{\frac{\sum_{{\bf x} \in X} (u({\bf x}) - U({\bf x};\theta^*))^2}{\sum_{{\bf x} \in X}  u({\bf x})^2} },
\end{equation}
where $U({\bf x};\theta^*) = U_i ({\bf x};\theta_i^*)$ for ${\bf x} \in \Omega_i $ and
the test sample set $X$ is obtained from the evenly spaced 251,001 ($501 \times 501$) samples in the square domain $\Omega$.

\subsection{Convergence test with a smooth solution}
In this subsection, we test and compare the convergence of the proposed Algorithms~\ref{algo1}-\ref{algo3}.
We consider the Poisson model problem~\eqref{model:Poisson} in a unit square domain $\Omega = (0,1)^2 $,
where the functions $f$ and $g$ are chosen to give a smooth exact solution,
\begin{equation}\label{smooth_test}
u({\bf x}):=u(x,y) = \sin(2\pi x) \sin (2\pi y).
\end{equation}
Before we proceed, we solve the above model problem by introducing a single neural network function $U({\bf x};\theta)$ and training its parameters $\theta$ for the standard PINN loss function and the one with the augmented Lagrangian term included for the boundary constraints, $u({\mathbf{x}})=g({\mathbf{x}})$,
in order to show the effectiveness of the augmented Lagrangian method
for enhancing the training performance and solution accuracy.

\begin{figure}[htb!]
 \centering
\includegraphics[height=4cm]{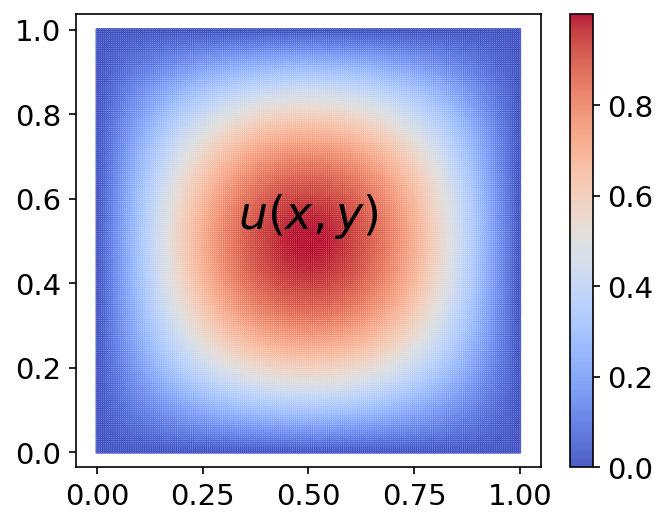}
\includegraphics[height=4cm]{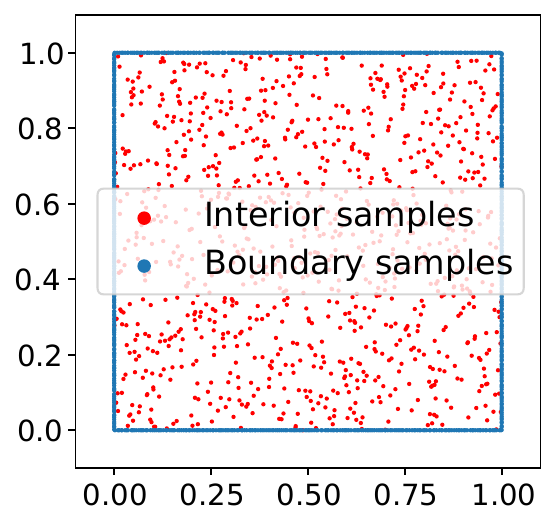}
  \caption{Plot of the exact solution~\eqref{smooth_test} (left)
  and an example of training samples for the single neural network (right).}
  \label{sec31_fig1}
  \end{figure}

The network and training data settings in the single neural network function are given as follows:
\begin{itemize}
\item {\bf Network structure} : Fully connected neural network of 4 hidden layers with 50 neurons per hidden layer.
\item {\bf The number of total parameters (weights, bias) of Network} : 7,851.
\item {\bf Training samples} : 1,000 samples for interior and 800 samples for boundary (200 for each edge).
\end{itemize}
The exact solution plot and an example of training sample points are presented in Figure~\ref{sec31_fig1}.

The two loss functions are defined as $\mathcal{J}_1 (\theta) :=  \mathcal{L}_{f} + \mathcal{L}_{g}$ and $\mathcal{J}_2 (\theta,\lambda_0) :=  \mathcal{L}_{f} + \mathcal{L}_{g} + \Lambda_{0}$,
where $\mathcal{L}_{f}$ is the PDE residual loss, $\mathcal{L}_{g}$ is the boundary loss, and $\Lambda_{0}$ is the augmented Lagrangian term for the boundary condition constraints,
$$\Lambda_0=\sum_{\mathbf{x} \in X_{\partial \Omega}} (U(\mathbf{x};\theta)-g(\mathbf{x})) \lambda_0(\mathbf{x}).$$
When optimizing $\lambda_0$ for $\mathcal{J}_2 (\theta,\lambda_0)$, we update $\lambda_0$ by the gradient ascent method every training epoch $\ell$,
\begin{equation}
\lambda_{0}^{(\ell)} (\mathbf{x}) = \lambda_{0}^{(\ell-1)} (\mathbf{x}) + \alpha_0 ( U (\mathbf{x}; \theta^{(\ell)}) - g(\mathbf{x}) )
\end{equation}
with $\alpha_0 = 0.1$.
In our computation, we observed that the use of a larger learning rate $\alpha_0$ than $lr(=0.001)$ gives a better training result for $\lambda_0(\mathbf{x})$ and we thus choose such a larger value of $\alpha_0$.

In Table~\ref{sec31_tb1}, the relative $L^2$-error values $\epsilon_u$ are reported for the two loss function cases.
The mean value of the errors in neural network approximate solutions for five different training sample sets
is reported and the standard deviation of the five error values is reported inside the parenthesis.
The errors are reduced as performing more training epochs in both cases.
The loss function $\mathcal{J}_2$ enriched with the augmented Lagrangian term
gives more accurate neural network approximate solutions than the simple PINN loss case, $\mathcal{J}_1$.
In the following numerical experiments, we will thus enrich
the localized loss function with the augmented Lagrangian term
in Algorithm~\ref{algo2} and Algorithm~\ref{algo3}
to increase the partitioned neural network solution accuracy further.

\begin{table}[h!]
\begin{center}
  \caption{Error values ${\epsilon}_u$ of the Poisson problem~\eqref{model:Poisson} with the smooth solution~\eqref{smooth_test}: Training results of the single neural network function for the standard PINN loss $\mathcal{J}_1(\theta)$
  and for the loss $\mathcal{J}_2(\theta,\lambda_0)$ enriched with the augmented Lagrangian term.
  The numbers are the mean value of errors (and standard deviation) obtained from five different training sample sets.
 }
{\footnotesize \renewcommand{\arraystretch}{1.2}
    \begin{tabular}{cccccc}
        \Xhline{3\arrayrulewidth}
Loss  & 10,000 epochs & 20,000 epochs & 50,000 epochs & 100,000 epochs
\\[1mm]         \Xhline{0.8\arrayrulewidth}
$\mathcal{J}_1 (\theta)$&
0.00577  {\scriptsize(1.39e-03)}&
0.00418  {\scriptsize(5.35e-04)}&
0.00313  {\scriptsize(5.96e-04)}&
0.00247  {\scriptsize(4.75e-04)} \\
$\mathcal{J}_2 (\theta,\lambda_0)$ &
0.00129  {\scriptsize(3.17e-04)}&
0.00070  {\scriptsize(2.16e-04)}&
0.00038  {\scriptsize(9.26e-05)}&
0.00037  {\scriptsize(1.19e-04)} \\
[0.5mm]     \Xhline{3\arrayrulewidth}
    \end{tabular}
}
\label{sec31_tb1}
\vskip-.7truecm
\end{center}
\end{table}

We now partition the domain $\Omega$ into non-overlapping subdomains $\Omega_i$, such that
$$ \overline{\Omega} = \overline{\cup_i \Omega_i}, \quad \Omega_i \cap \Omega_j = \emptyset, \quad \Gamma_{ij} = \partial \Omega_i \cap \partial \Omega_j , \quad \Gamma = \cup_{ij} \Gamma_{ij}. $$
Here $\Gamma_{ij}$ denotes the interface of two subdomains and $\Gamma$ denotes the interface in the subdomain
partition. The model problem~\eqref{model:Poisson} is then partitioned into the local problems~\eqref{eq_poi_dd},
where we note that the interface jump condition are set to the zero value for the smooth solution~\eqref{smooth_test},
i.e., $p_{ij}=q_{ij}=0$.

\begin{figure}[htb!]
 \centering
 \includegraphics[height=3.7cm]{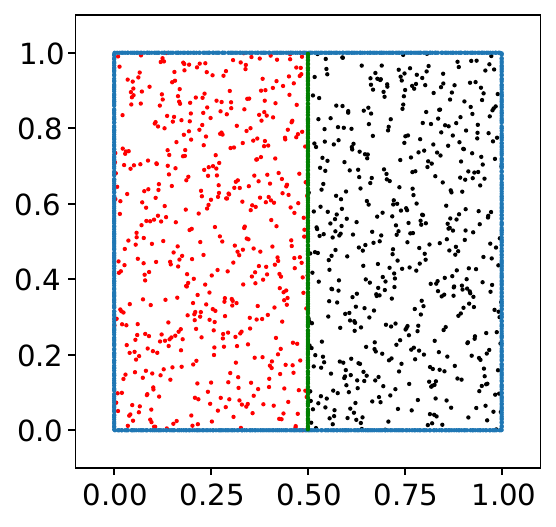}
  \includegraphics[height=3.7cm]{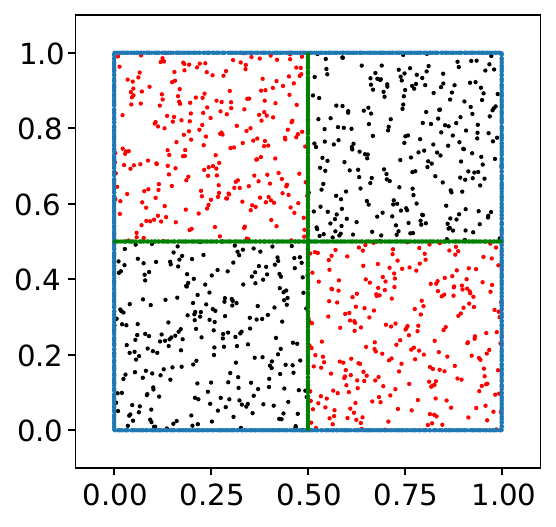}
  \includegraphics[height=3.7cm]{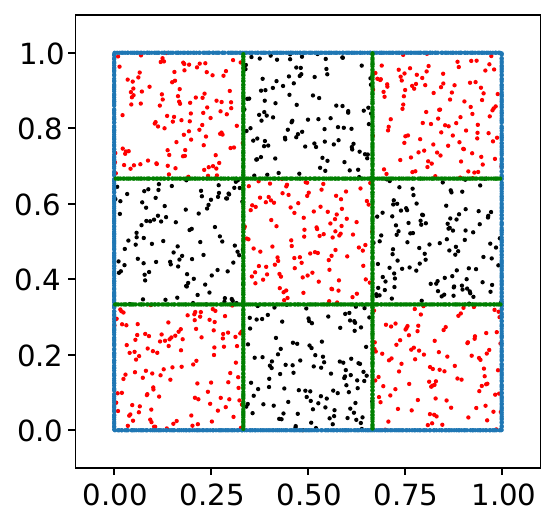}
  \includegraphics[height=3.7cm]{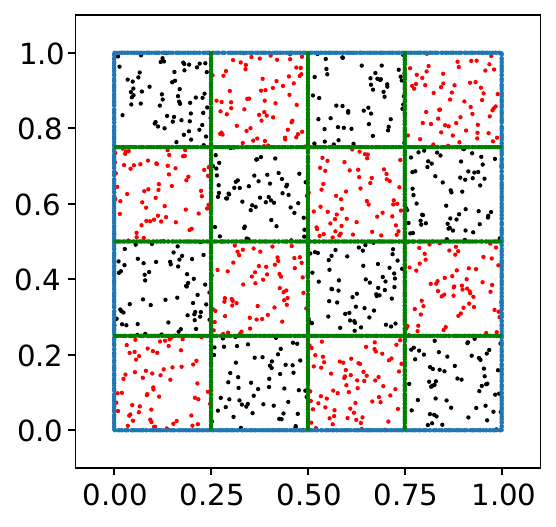}
  \caption{Plots of training samples and subdomain partitions: 2, 4, 9, and 16 subdomain partitions (left to right).}
  \label{sec31_fig2}
  \end{figure}

As shown in Figure~\ref{sec31_fig2}, we consider 2, 4, 9, and 16 subdomain partitions
in our numerical experiments.
For all the subdomain partition case, we introduce the local neural network as a fully connected neural network of
four hidden layers and the number of neurons per hidden layer as described in Table~\ref{network_tbl}.
\begin{table}[h!]
\begin{center}
  \caption{Information of the local neural networks for all the subdomain partition case: The number of neurons per each hidden layer and the number of parameters per each local neural network. }
{\footnotesize \renewcommand{\arraystretch}{1.2}
    \begin{tabular}{ccc}
        \Xhline{3\arrayrulewidth}
Subdomains  & Neurons & Parameters
\\[1mm]         \Xhline{0.8\arrayrulewidth}
2 & 35 & 3,921 \\
4 & 23 & 1,749 \\
9 & 16 & 881 \\
16 & 11 & 488 \\
[0.5mm]     \Xhline{3\arrayrulewidth}
    \end{tabular}
}
\label{network_tbl}
\end{center}
\end{table}

In addition, we set 2,000 training samples in the entire domain $\Omega$ using a Latin hypercube sampling method,
and we then choose each subdomain training samples from those 2,000 samples that belong to the subdomain. For the interface sampling points, 200 training samples are introduced for each vertical and horizontal interfaces, resulting in 200, 400, 800, and 1,200 interface samples totally for 2, 4, 9, and 16 subdomain partition cases, respectively.
An example of training samples for the given subdomain partition
is also presented for $2$, $4$, $9$, and $16$ subdomain partition cases
in Figure~\ref{sec31_fig2}.

In Table~\ref{sec31_tb2}, we report the test errors of the training results of the three algorithms
for 2, 4, 9, and 16 subdomain partition cases, respectively.
The mean value and the standard deviation of training results from five different training samples are reported.
The errors are observed at $10,000$, $20,000$, $50,000$, and $100,000$ total training epochs to compare the training performance
of the three algorithms.
In Algorithm~\ref{algo2}, the augmented Lagrangian term is included in the loss function.
The error results show that such an enriched loss function gives better training results than
those obtained from Algorithm~\ref{algo1}.
As more subdomains are introduced in the partition, i.e., the subdomain interface $\Gamma$ becomes larger,
the augmented Lagrangian term works more effectively to give smaller errors than in Algorithm~\ref{algo1}.
In Algorithm~\ref{algo2}, the data communication between neighboring neural network solutions is needed
every training epoch.
To reduce such a data communication cost, we proposed Algorithm~\ref{algo3}.
To test the convergence of Algorithm~\ref{algo3},
we perform Algorithm~\ref{algo3} with different settings for $N_l$.
We choose $N_l$ as $100$, $200$, $500$, and $1,000$ and report the error values depending on the choice of $N_l$.
With a larger $N_l$, a less number of outer iteration is needed and the data communication cost is thus reduced.
However, the less communication case resulted in slightly larger errors than those obtained from more data communication cases at the early outer iterations.
As proceeding more outer iterations, up to  the $100,000$ training epochs, Algorithm~\ref{algo3} with $N_l=1,000$ gives still acceptable error results for all the subdomain partition cases.

\begin{table}[h!]
\begin{center}
  \caption{Error values ${\epsilon}_u$ of the Poisson problem~\eqref{model:Poisson} with the smooth solution~\eqref{smooth_test}: The learning rate $\alpha_0 = 0.1$ for algorithms $\mathcal{A}_2$ and  $\mathcal{A}_3$.
  $N_l$ is the number of training epochs of local problems in algorithm $\mathcal{A}_3$ for every outer iteration.
  The numbers are the mean value of errors (and standard deviation) (below: the number of data communication) obtained from five different training sample sets.
 } \label{sec31_tb2}
{\footnotesize \renewcommand{\arraystretch}{1.2}
    \begin{tabular}{clllll}
        \Xhline{3\arrayrulewidth}
\multicolumn{2}{c}{\bf \small Two subdomains}  & Max epochs &  &  &  \\
Algorithm & Parameters  & 10,000 epochs & 20,000 epochs & 50,000 epochs & 100,000 epochs
\\[1mm]         \Xhline{0.8\arrayrulewidth}
$\mathcal{A}_1$& &
0.01159  {\scriptsize(3.55e-03)} \vspace{-0.15cm}&
0.00649  {\scriptsize(2.17e-03)}&
0.00356  {\scriptsize(7.04e-04)}&
0.00240  {\scriptsize(8.46e-04)} \\
& & \vspace{-0.01cm} {\scriptsize(10000)}
& \vspace{-0.01cm} {\scriptsize(20000)}
& \vspace{-0.01cm} {\scriptsize(50000)}
& \vspace{-0.01cm} {\scriptsize(100000)}
\\ \Xhline{0.8\arrayrulewidth}
$\mathcal{A}_2$ & $\alpha_\lambda=0.1 $&
0.00165  {\scriptsize(8.64e-04)} \vspace{-0.15cm}&
0.00080 {\scriptsize(2.55e-04)}&
0.00040  {\scriptsize(3.88e-05)}&
0.00035  {\scriptsize(1.06e-04)} \\
& & \vspace{-0.01cm} {\scriptsize(10000)}
& \vspace{-0.01cm} {\scriptsize(20000)}
& \vspace{-0.01cm} {\scriptsize(50000)}
& \vspace{-0.01cm} {\scriptsize(100000)}
\\ \Xhline{0.8\arrayrulewidth}
$\mathcal{A}_3$ & $N_l = 100, \alpha_\lambda = 0.1$&
0.00341  {\scriptsize(9.29e-04)}  \vspace{-0.15cm}&
0.00203  {\scriptsize(9.82e-04)}&
0.00100  {\scriptsize(3.52e-04)}&
0.00100  {\scriptsize(1.87e-04)} \\
& & \vspace{-0.01cm} {\scriptsize(100)}
& \vspace{-0.01cm} {\scriptsize(200)}
& \vspace{-0.01cm} {\scriptsize(500)}
& \vspace{-0.01cm} {\scriptsize(1000)}
\\
$\mathcal{A}_3$ & $N_l = 200, \alpha_\lambda = 0.05$&
0.00219  {\scriptsize(6.70e-04)}  \vspace{-0.15cm}&
0.00100  {\scriptsize(1.12e-04)}&
0.00077  {\scriptsize(2.03e-04)}&
0.00077  {\scriptsize(4.82e-04)} \\
& & \vspace{-0.01cm} {\scriptsize(50)}
& \vspace{-0.01cm} {\scriptsize(100)}
& \vspace{-0.01cm} {\scriptsize(250)}
& \vspace{-0.01cm} {\scriptsize(500)}
\\
$\mathcal{A}_3$ & $N_l = 500, \alpha_\lambda = 0.01$&
0.00236  {\scriptsize(3.49e-04)}  \vspace{-0.15cm}&
0.00159  {\scriptsize(7.03e-04)}&
0.00091  {\scriptsize(2.28e-04)}&
0.00075  {\scriptsize(2.94e-04)} \\
& & \vspace{-0.01cm} {\scriptsize(20)}
& \vspace{-0.01cm} {\scriptsize(40)}
& \vspace{-0.01cm} {\scriptsize(100)}
& \vspace{-0.01cm} {\scriptsize(200)}
\\
$\mathcal{A}_3$ & $N_l = 1,000, \alpha_\lambda = 0.005$&
0.00385  {\scriptsize(1.74e-03)}  \vspace{-0.15cm}&
0.00146  {\scriptsize(5.30e-04)}&
0.00137  {\scriptsize(4.25e-04)}&
0.00066  {\scriptsize(2.20e-04)} \\
& & \vspace{-0.01cm} {\scriptsize(10)}
& \vspace{-0.01cm} {\scriptsize(20)}
& \vspace{-0.01cm} {\scriptsize(50)}
& \vspace{-0.01cm} {\scriptsize(100)}
    \end{tabular}
}
{\footnotesize \renewcommand{\arraystretch}{1.2}
    \begin{tabular}{clllll}
        \Xhline{3\arrayrulewidth}
\multicolumn{2}{c}{\bf \small Four subdomains}  & Max epochs &  &  &  \\
Algorithm & Parameters  & 10,000 epochs & 20,000 epochs & 50,000 epochs & 100,000 epochs
\\[1mm]         \Xhline{0.8\arrayrulewidth}
$\mathcal{A}_1$& &
0.01131  {\scriptsize(3.06e-03)} \vspace{-0.15cm}&
0.00619  {\scriptsize(2.08e-03)}&
0.00352  {\scriptsize(8.88e-04)}&
0.00224  {\scriptsize(2.31e-04)} \\
& & \vspace{-0.01cm} {\scriptsize(10000)}
& \vspace{-0.01cm} {\scriptsize(20000)}
& \vspace{-0.01cm} {\scriptsize(50000)}
& \vspace{-0.01cm} {\scriptsize(100000)}
\\ \Xhline{0.8\arrayrulewidth}
$\mathcal{A}_2$& $\alpha_\lambda = 0.1 $&
0.00278  {\scriptsize(1.27e-03)} \vspace{-0.15cm}&
0.00089  {\scriptsize(3.22e-04)}&
0.00059  {\scriptsize(4.03e-04)}&
0.00030  {\scriptsize(5.71e-05)} \\
& & \vspace{-0.01cm} {\scriptsize(10000)}
& \vspace{-0.01cm} {\scriptsize(20000)}
& \vspace{-0.01cm} {\scriptsize(50000)}
& \vspace{-0.01cm} {\scriptsize(100000)}
\\ \Xhline{0.8\arrayrulewidth}
$\mathcal{A}_3$& $N_l = 100, \alpha_\lambda = 0.1$&
0.00470  {\scriptsize(1.06e-03)}  \vspace{-0.15cm}&
0.00184  {\scriptsize(6.19e-04)}&
0.00146  {\scriptsize(5.38e-04)}&
0.00121  {\scriptsize(4.49e-04)} \\
& & \vspace{-0.01cm} {\scriptsize(100)}
& \vspace{-0.01cm} {\scriptsize(200)}
& \vspace{-0.01cm} {\scriptsize(500)}
& \vspace{-0.01cm} {\scriptsize(1000)}
\\
$\mathcal{A}_3$& $N_l = 200, \alpha_\lambda = 0.05$&
0.00512  {\scriptsize(1.50e-03)}  \vspace{-0.15cm}&
0.00167  {\scriptsize(5.20e-04)}&
0.00134  {\scriptsize(6.26e-04)}&
0.00094  {\scriptsize(2.29e-04)} \\
& & \vspace{-0.01cm} {\scriptsize(50)}
& \vspace{-0.01cm} {\scriptsize(100)}
& \vspace{-0.01cm} {\scriptsize(200)}
& \vspace{-0.01cm} {\scriptsize(500)}
\\
$\mathcal{A}_3$& $N_l = 500, \alpha_\lambda = 0.01$&
0.00713  {\scriptsize(2.15e-03)}  \vspace{-0.15cm}&
0.00346  {\scriptsize(1.49e-03)}&
0.00165  {\scriptsize(9.51e-04)}&
0.00139  {\scriptsize(8.76e-04)} \\
& & \vspace{-0.01cm} {\scriptsize(20)}
& \vspace{-0.01cm} {\scriptsize(40)}
& \vspace{-0.01cm} {\scriptsize(100)}
& \vspace{-0.01cm} {\scriptsize(200)}
\\
$\mathcal{A}_3$& $N_l = 1,000, \alpha_\lambda = 0.005$&
0.01093 {\scriptsize(3.55e-03)}  \vspace{-0.15cm}&
0.00557  {\scriptsize(1.75e-03)}&
0.00259  {\scriptsize(7.08e-04)}&
0.00154 {\scriptsize(5.30e-04)} \\
& & \vspace{-0.01cm} {\scriptsize(10)}
& \vspace{-0.01cm} {\scriptsize(20)}
& \vspace{-0.01cm} {\scriptsize(50)}
& \vspace{-0.01cm} {\scriptsize(100)}
    \end{tabular}
}
{\footnotesize \renewcommand{\arraystretch}{1.2}
    \begin{tabular}{clllll}
        \Xhline{3\arrayrulewidth}
\multicolumn{2}{c}{\bf \small Nine subdomains}  & Max epochs &  &  &  \\
Algorithm & Parameters  & 10,000 epochs & 20,000 epochs & 50,000 epochs & 100,000 epochs
\\[1mm]         \Xhline{0.8\arrayrulewidth}
$\mathcal{A}_1$& &
0.02468  {\scriptsize(2.71e-03)} \vspace{-0.15cm}&
0.01169  {\scriptsize(2.10e-03)}&
0.00511  {\scriptsize(7.14e-04)}&
0.00386  {\scriptsize(2.02e-04)} \\
& & \vspace{-0.01cm} {\scriptsize(10000)}
& \vspace{-0.01cm} {\scriptsize(20000)}
& \vspace{-0.01cm} {\scriptsize(50000)}
& \vspace{-0.01cm} {\scriptsize(100000)}
\\ \Xhline{0.8\arrayrulewidth}
$\mathcal{A}_2$& $\alpha_\lambda=0.1 $&
0.00330  {\scriptsize(6.69e-04)} \vspace{-0.15cm}&
0.00175  {\scriptsize(5.58e-04)}&
0.00071  {\scriptsize(3.64e-04)}&
0.00041  {\scriptsize(1.00e-04)} \\
& & \vspace{-0.01cm} {\scriptsize(10000)}
& \vspace{-0.01cm} {\scriptsize(20000)}
& \vspace{-0.01cm} {\scriptsize(50000)}
& \vspace{-0.01cm} {\scriptsize(100000)}
\\ \Xhline{0.8\arrayrulewidth}
$\mathcal{A}_3$& $N_l = 100, \alpha_\lambda = 0.1$&
0.00815  {\scriptsize(1.48e-03)}  \vspace{-0.15cm}&
0.00257  {\scriptsize(2.99e-04)}&
0.00174  {\scriptsize(2.68e-04)}&
0.00110  {\scriptsize(2.64e-04)} \\
& & \vspace{-0.01cm} {\scriptsize(100)}
& \vspace{-0.01cm} {\scriptsize(200)}
& \vspace{-0.01cm} {\scriptsize(500)}
& \vspace{-0.01cm} {\scriptsize(1000)}
\\
$\mathcal{A}_3$& $N_l = 200, \alpha_\lambda = 0.05$&
0.00982  {\scriptsize(3.21e-03)}  \vspace{-0.15cm}&
0.00297  {\scriptsize(3.39e-04)}&
0.00144  {\scriptsize(5.54e-04)}&
0.00131  {\scriptsize(1.81e-04)} \\
& & \vspace{-0.01cm} {\scriptsize(50)}
& \vspace{-0.01cm} {\scriptsize(100)}
& \vspace{-0.01cm} {\scriptsize(200)}
& \vspace{-0.01cm} {\scriptsize(500)}
\\
$\mathcal{A}_3$& $N_l = 500, \alpha_\lambda = 0.01$&
0.01551  {\scriptsize(3.12e-03)}  \vspace{-0.15cm}&
0.00487  {\scriptsize(1.47e-03)}&
0.00206  {\scriptsize(2.57e-04)}&
0.00117  {\scriptsize(3.51e-05)} \\
& & \vspace{-0.01cm} {\scriptsize(20)}
& \vspace{-0.01cm} {\scriptsize(40)}
& \vspace{-0.01cm} {\scriptsize(100)}
& \vspace{-0.01cm} {\scriptsize(200)}
\\
$\mathcal{A}_3$& $N_l = 1,000, \alpha_\lambda = 0.005$&
0.09110  {\scriptsize(1.41e-02)}  \vspace{-0.15cm}&
0.01193  {\scriptsize(1.58e-03)}&
0.00314  {\scriptsize(3.22e-04)}&
0.00148  {\scriptsize(1.53e-04)} \\
& & \vspace{-0.01cm} {\scriptsize(10)}
& \vspace{-0.01cm} {\scriptsize(20)}
& \vspace{-0.01cm} {\scriptsize(50)}
& \vspace{-0.01cm} {\scriptsize(100)}
    \end{tabular}
}
{\footnotesize \renewcommand{\arraystretch}{1.2}
    \begin{tabular}{clllll}
        \Xhline{3\arrayrulewidth}
\multicolumn{2}{c}{\bf \small Sixteen subdomains}  & Max epochs &  &  &  \\
Algorithm & Parameters  & 10,000 epochs & 20,000 epochs & 50,000 epochs & 100,000 epochs
\\[1mm]         \Xhline{0.8\arrayrulewidth}
$\mathcal{A}_1$& &
0.02902  {\scriptsize(4.49e-03)} \vspace{-0.15cm}&
0.01764  {\scriptsize(2.48e-03)}&
0.00885  {\scriptsize(1.14e-03)}&
0.00551  {\scriptsize(6.43e-04)} \\
& & \vspace{-0.01cm} {\scriptsize(10000)}
& \vspace{-0.01cm} {\scriptsize(20000)}
& \vspace{-0.01cm} {\scriptsize(50000)}
& \vspace{-0.01cm} {\scriptsize(100000)}
\\ \Xhline{0.8\arrayrulewidth}
$\mathcal{A}_2$& $\alpha_\lambda=0.1 $&
0.00373 {\scriptsize(6.87e-04)} \vspace{-0.15cm}&
0.00152  {\scriptsize(1.38e-04)}&
0.00078  {\scriptsize(2.37e-04)}&
0.00056  {\scriptsize(2.20e-04)} \\
& & \vspace{-0.01cm} {\scriptsize(10000)}
& \vspace{-0.01cm} {\scriptsize(20000)}
& \vspace{-0.01cm} {\scriptsize(50000)}
& \vspace{-0.01cm} {\scriptsize(100000)}
\\ \Xhline{0.8\arrayrulewidth}
$\mathcal{A}_3$& $N_l = 100, \alpha_\lambda = 0.1$&
0.00834  {\scriptsize(9.78e-04)}  \vspace{-0.15cm}&
0.00349  {\scriptsize(5.13e-04)}&
0.00240  {\scriptsize(7.69e-04)}&
0.00155  {\scriptsize(2.45e-04)} \\
& & \vspace{-0.01cm} {\scriptsize(100)}
& \vspace{-0.01cm} {\scriptsize(200)}
& \vspace{-0.01cm} {\scriptsize(500)}
& \vspace{-0.01cm} {\scriptsize(1000)}
\\
$\mathcal{A}_3$& $N_l = 200, \alpha_\lambda = 0.05$&
0.02477  {\scriptsize(3.72e-03)}  \vspace{-0.15cm}&
0.00430  {\scriptsize(7.77e-04)}&
0.00213  {\scriptsize(1.40e-04)}&
0.00177  {\scriptsize(4.79e-04)} \\
& & \vspace{-0.01cm} {\scriptsize(50)}
& \vspace{-0.01cm} {\scriptsize(100)}
& \vspace{-0.01cm} {\scriptsize(200)}
& \vspace{-0.01cm} {\scriptsize(500)}
\\
$\mathcal{A}_3$& $N_l = 500, \alpha_\lambda = 0.01$&
0.11664  {\scriptsize(1.38e-02)}  \vspace{-0.15cm}&
0.01224  {\scriptsize(1.35e-03)}&
0.00294  {\scriptsize(6.02e-04)}&
0.00197  {\scriptsize(3.87e-04)} \\
& & \vspace{-0.01cm} {\scriptsize(20)}
& \vspace{-0.01cm} {\scriptsize(40)}
& \vspace{-0.01cm} {\scriptsize(100)}
& \vspace{-0.01cm} {\scriptsize(200)}
\\
$\mathcal{A}_3$& $N_l = 1,000, \alpha_\lambda = 0.005$&
0.34196  {\scriptsize(6.59e-02)}  \vspace{-0.15cm}&
0.10457   {\scriptsize(2.02e-02)}&
0.00514   {\scriptsize(3.32e-04)}&
0.00226   {\scriptsize(5.19e-04)} \\
& & \vspace{-0.01cm} {\scriptsize(10)}
& \vspace{-0.01cm} {\scriptsize(20)}
& \vspace{-0.01cm} {\scriptsize(50)}
& \vspace{-0.01cm} {\scriptsize(100)}
\\
[0.5mm]     \Xhline{3\arrayrulewidth}
    \end{tabular}
}
\end{center}
\end{table}

\subsection{Elliptic problem with discontinuous coefficient}
\begin{figure}[htb!]
 \centering
  \includegraphics[height=4cm]{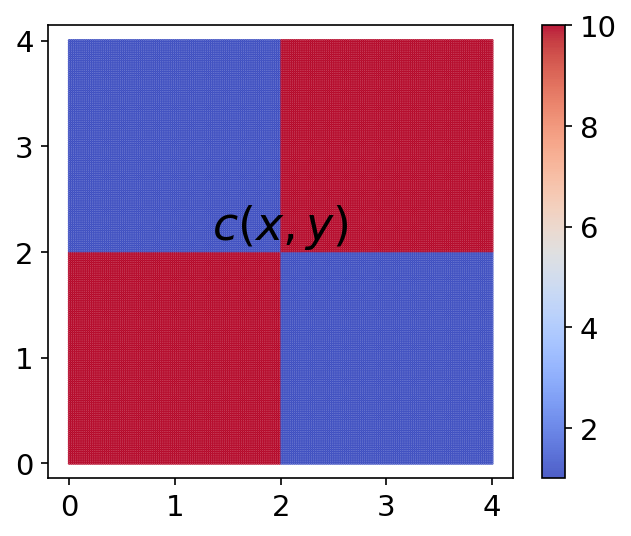}
  \includegraphics[height=4cm]{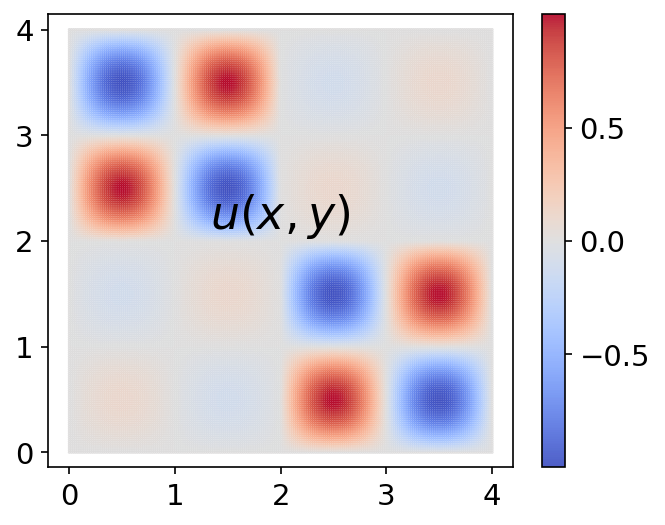}
  \includegraphics[height=4cm]{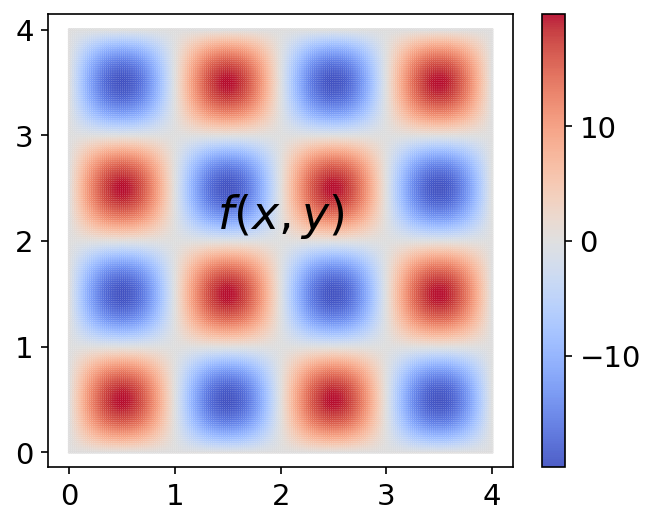}
  \caption{Coefficient (left), exact solution (middle) and force term (right) of 
  the model problem~\eqref{eq_disc_coef}-\eqref{disc_coef}.}
  \label{sec32_fig1}
  \end{figure}

In this subsection, we consider a more challenging test example, i.e., a second-order elliptic problem with discontinuous coefficients.
For the computational domain $\Omega = (0,4)^2$, the governing equation is given as follows:
\begin{equation}\label{eq_disc_coef}
\begin{split}
-\nabla \cdot (c(\mathbf{x}) \nabla u(\mathbf{x})) &=f(\mathbf{x}) \; \text{ in } \Omega, \\
u(\mathbf{x}) &= 0 \; \text{ on } \partial\Omega,
\end{split}
\end{equation}
where $f({\bf x}):=f(x,y) = 2 \pi^2 \sin( \pi x) \sin (\pi y)$ and the coefficient $c({\bf x}):=c(x,y)$ is piecewise constant with respect to the partition of the domain \begin{equation}\label{disc_coef}
c(x,y) = \left\{
\begin{array}{cc}
10 &  \text{in } \Omega_1 \cup \Omega_4, \\
1 &  \text{in } \Omega_2 \cup \Omega_3,
\end{array}
\right.
\end{equation}
where $\Omega_1 = (0,2)^2$, $\Omega_2 = (2,4) \times (0,2)$, $\Omega_3 = (0,2) \times (2,4)$ and , $\Omega_4 = (2,4)^2$. The exact solution is then obtained as $u(x,y) = \sin(\pi x) \sin (\pi y) / c(x,y)$.
The plots of the discontinuous coefficient value, the exact solution, and the right hand side forcing term $f$
are presented in Figure~\ref{sec32_fig1}.

\begin{figure}[htb!]
 \centering
  \includegraphics[height=4cm]{figs/subsec32/01_exact_u.png}
  \includegraphics[height=4cm]{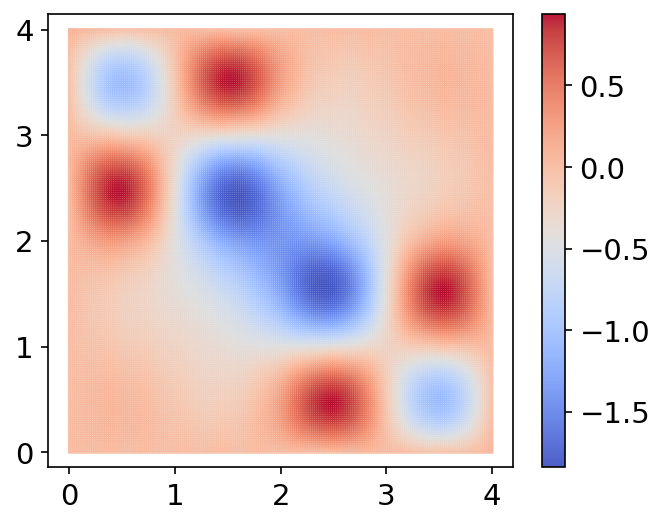}
  \includegraphics[height=4cm]{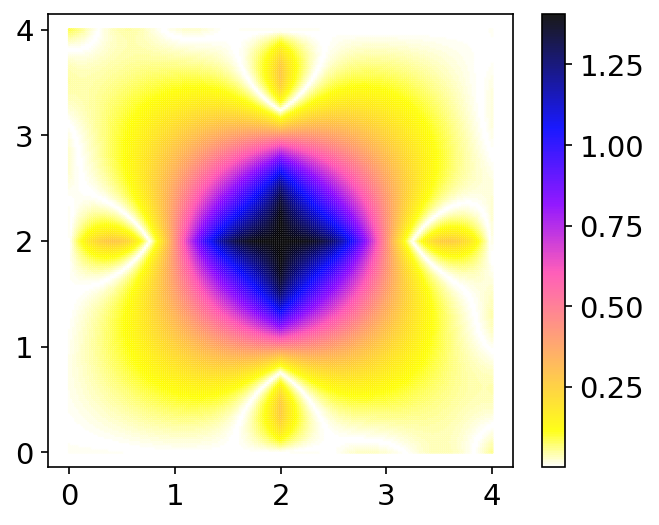}
  \caption{Target solution (left), trained solution with a single neural network (middle), and point-wise errors (right) in the trained solution.}
  \label{sec32_fig2}
  \end{figure}

Such a discontinuous coefficient model problem is difficult to be solved by a single neural network. In order to investigate the feasibility of solving the given problem~\eqref{eq_disc_coef}--\eqref{disc_coef} with a single neural network, we employed the following network structure and hyperparameter settings to approximate the model solution.
\begin{itemize}
\item {\bf Network structure} : Fully connected neural network of 4 hidden layers with 50 neurons per hidden layer.
\item {\bf The number of total parameters (weights, bias) of Network} : 7,851.
\item {\bf Training samples} : 2,000 samples for interior and 800 samples for boundary (200 for each edge).
\end{itemize}
For the parameter optimization, we used the standard PINN loss function with the PDE residual loss and boundary condition loss terms.
In Figure~\ref{sec32_fig2}, the error plot shows that the single neural network approximation gives larger errors near the interface of the coefficient discontinuity. The single neural network is thus not suitable for handling such discontinuous coefficient model problems.

For the same discontinuous coefficient model problem, we partition the domain into four non-overlapping subdomains
so that the subdomain interface can handle the model coefficient discontinuity easily.
The above single domain problem~\eqref{eq_disc_coef}-\eqref{disc_coef} is then partitioned into
the following subdomain problems:
\begin{equation}\label{eq_disc_coef_if}
\begin{split}
-\nabla \cdot (c(\mathbf{x}) \nabla u_i(\mathbf{x})) &=f(\mathbf{x}) \; \text{ in } \Omega_i, \\
u_i(\mathbf{x}) &= 0 \; \text{ on } \partial\Omega_i, \\
\left\llbracket  u_i ({\bf x})\right\rrbracket_{\Gamma_{ij}} &= 0 \; \text{ on } \Gamma_{ij},  \ \ \forall j \in s(i), \\
\left\llbracket  c({\bf x}) \frac{\partial u_i}{\partial n}({\bf x}) \right\rrbracket _{\Gamma_{ij}} &= 0 \; \text{ on } \Gamma_{ij},  \ \ \forall j \in s(i),
\end{split}
\end{equation}
where we note that the coefficient $c(\mathbf{x})$ does not have discontinuity inside $\Omega_i$,
and the coefficient discontinuity only appears in the evaluation of the normal flux condition on the interface $\Gamma_{ij}$.
For each subdomain, we employ the corresponding local neural network function to approximate the model solution.
The neural network configuration and training sample information are as follows:
\begin{itemize}
\item {\bf Network structure} : 4 hidden layers with 23 neurons per hidden layer in each local neural network.
\item {\bf The number of total parameters (weights, bias) of Network} : 1,749 in each local neural network.
\item {\bf Training samples in $\Omega$} : 2,000 samples for interior, 800 samples for boundary (200 for each boundary edge) and 400 samples for interface (200 for the vertical interface line and 200 for the horizontal interface line).
\end{itemize}

The error results are presented in Table~\ref{sec32_tb1}.
Algorithm~\ref{algo1} without the augmented Lagrangian term is observed to give
larger errors than those in Algorithms~\ref{algo2} and \ref{algo3}.
In Algorithm~\ref{algo2}, we obtained the best results where the data exchange between
the neighboring neural networks is needed every training epoch.
The error results with various choices of $N_l$ in Algorithm~\ref{algo3}, i.e., less data communication than in Algorithm~\ref{algo2}, show similar results as in the previous smooth example. At the $100,000$ total epochs,
with a much lesser communication of $100$ outer iterations with $N_l=1,000$, compared to $1,000$ outer iterations with $N_l=100$, we obtained still acceptable error results, even smaller than that in Algorithm~\ref{algo1}.

\begin{table}[htb!]
\begin{center}
  \caption{Error values $\epsilon_u$ of the discontinuous coefficient elliptic problem in \eqref{eq_disc_coef}-\eqref{disc_coef} with four local neural networks:
   The learning rate $\alpha_0 = 0.1$ is used for algorithms $\mathcal{A}_2$ and  $\mathcal{A}_3$.
   $N_l$ is the number of training epochs of local problems in algorithm~$\mathcal{A}_3$ for every outer iteration.
   The numbers are the mean value of errors (and standard deviation) (below: the number of data communication) obtained from five different training sample sets.
  }
{\footnotesize \renewcommand{\arraystretch}{1.2}
    \begin{tabular}{clllll}
        \Xhline{3\arrayrulewidth}
\multicolumn{2}{c}{\bf Discontinuous coefficient}  & Max epochs &  &  &  \\
Algorithm & Parameters  & 10,000 epochs & 20,000 epochs & 50,000 epochs & 100,000 epochs
\\[1mm]         \Xhline{0.8\arrayrulewidth}
$\mathcal{A}_1$& &
0.05805  {\scriptsize(2.55e-02)} \vspace{-0.15cm}&
0.03308  {\scriptsize(1.03e-02)}&
0.01645  {\scriptsize(5.84e-03)}&
0.01186  {\scriptsize(4.20e-03)} \\
& & \vspace{-0.01cm} {\scriptsize(10000)}
& \vspace{-0.01cm} {\scriptsize(20000)}
& \vspace{-0.01cm} {\scriptsize(50000)}
& \vspace{-0.01cm} {\scriptsize(100000)}
\\ \Xhline{0.8\arrayrulewidth}
$\mathcal{A}_2$& $\alpha_\lambda=0.1 $&
0.00584 {\scriptsize(2.31e-03)} \vspace{-0.15cm}&
0.00373 {\scriptsize(1.74e-03)}&
0.00282 {\scriptsize(1.08e-03)}&
0.00123 {\scriptsize(9.43e-04)} \\
& & \vspace{-0.01cm} {\scriptsize(10000)}
& \vspace{-0.01cm} {\scriptsize(20000)}
& \vspace{-0.01cm} {\scriptsize(50000)}
& \vspace{-0.01cm} {\scriptsize(100000)}
\\ \Xhline{0.8\arrayrulewidth}
$\mathcal{A}_3$& $N_l = 100, \alpha_\lambda = 0.05$&
0.03732 {\scriptsize(1.56e-02)}  \vspace{-0.15cm}&
0.02484 {\scriptsize(9.64e-03)}&
0.01431 {\scriptsize(5.45e-03)}&
0.00795 {\scriptsize(3.46e-03)} \\
& & \vspace{-0.01cm} {\scriptsize(100)}
& \vspace{-0.01cm} {\scriptsize(200)}
& \vspace{-0.01cm} {\scriptsize(500)}
& \vspace{-0.01cm} {\scriptsize(1000)}
\\
$\mathcal{A}_3$& $N_l = 200, \alpha_\lambda = 0.01$&
0.04563  {\scriptsize(1.50e-02)}  \vspace{-0.15cm}&
0.02784  {\scriptsize(6.05e-03)}&
0.01699  {\scriptsize(3.31e-03)}&
0.00806  {\scriptsize(1.13e-03)} \\
& & \vspace{-0.01cm} {\scriptsize(50)}
& \vspace{-0.01cm} {\scriptsize(100)}
& \vspace{-0.01cm} {\scriptsize(200)}
& \vspace{-0.01cm} {\scriptsize(500)}
\\
$\mathcal{A}_3$& $N_l = 500, \alpha_\lambda = 0.002$&
0.04939   {\scriptsize(1.55e-02)}  \vspace{-0.15cm}&
0.03176   {\scriptsize(8.91e-03)}&
0.01893   {\scriptsize(6.28e-03)}&
0.00950   {\scriptsize(3.32e-03)} \\
& & \vspace{-0.01cm} {\scriptsize(20)}
& \vspace{-0.01cm} {\scriptsize(40)}
& \vspace{-0.01cm} {\scriptsize(100)}
& \vspace{-0.01cm} {\scriptsize(200)}
\\
$\mathcal{A}_3$& $N_l = 1,000, \alpha_\lambda = 0.001$&
0.04745  {\scriptsize(1.46e-02)}  \vspace{-0.15cm}&
0.03102   {\scriptsize(9.30e-03)}&
0.01807   {\scriptsize(5.21e-03)}&
0.00928  {\scriptsize(4.17e-03)} \\
& & \vspace{-0.01cm} {\scriptsize(10)}
& \vspace{-0.01cm} {\scriptsize(20)}
& \vspace{-0.01cm} {\scriptsize(50)}
& \vspace{-0.01cm} {\scriptsize(100)}
\\
[0.5mm]     \Xhline{3\arrayrulewidth}
    \end{tabular}
}
\label{sec32_tb1}
\end{center}
\end{table}

\subsection{Poisson interface problem}
\begin{figure}[htb!]
 \centering
  \includegraphics[height=4cm]{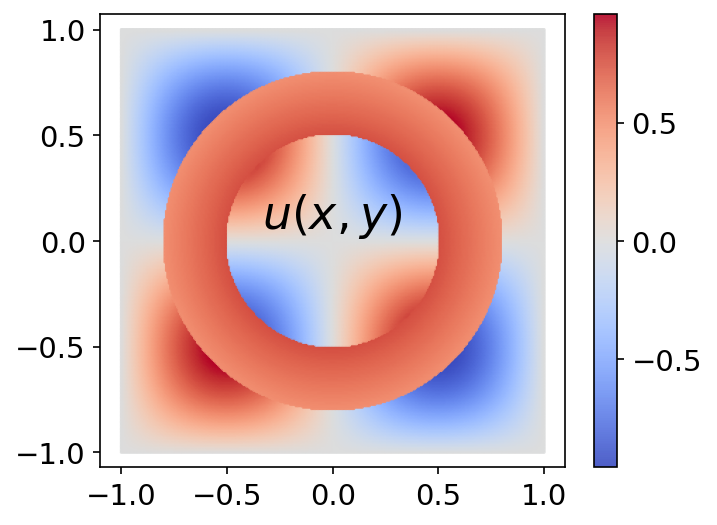}
  \includegraphics[height=4cm]{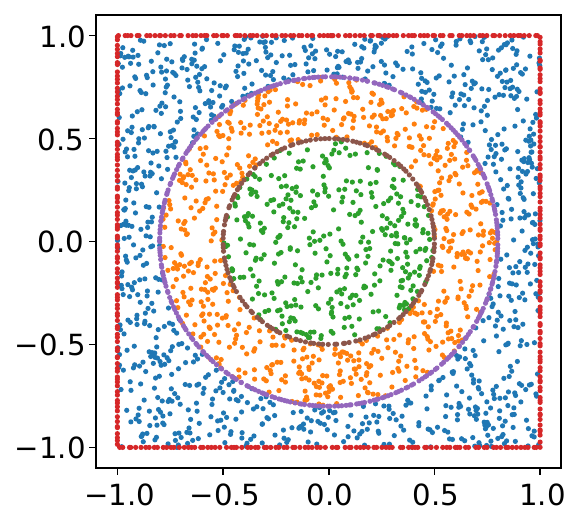}
  \caption{Exact solution~\eqref{ex_jump} (left) and
  an example of training sample points over the subdomain partition (right).}
  \label{sec33_fig1}
  \end{figure}

In this subsection, our goal is to solve the following Poisson interface problem subject to a Dirichlet boundary condition with non-vanishing jump values of $u$ on the interface $\Gamma$:
\begin{equation}\label{eq_poi_jump}
\begin{split}
-\triangle u ( {\bf x} ) & = f({\bf x})  \; \text{ in } \Omega, \\
u({\bf x}) &= g( {\bf x})  \; \text{ on } \partial \Omega, \\
\left\llbracket  u ({\bf x}) \right\rrbracket_{\Gamma} &= p( {\bf x})  \; \text{ on } \Gamma, \\
\left\llbracket  \frac{\partial u}{\partial { n}} ({\bf x}) \right\rrbracket _{\Gamma} &= q( {\bf x})  \; \text{ on } \Gamma,
\end{split}
\end{equation}
where the unit square domain $\Omega=(-1, \, 1)^2$ is partitioned into three non-overlapping subdomains, $\Omega_1 = \{ (x,y)\in \Omega \ : \ x^2+y^2<0.5^2\}$, $\Omega_2 = \{ (x,y)\in \Omega \ : \ 0.5^2 < x^2+y^2<0.8^2\}$, and $\Omega_3 = \{ (x,y)\in \Omega \ : \ 0.8^2 < x^2+y^2\}$, and their interface $\Gamma$ consists of $\Gamma_1$ and $\Gamma_2$, that are circles
centered at the origin with radius $r=0.5$ for $\Gamma_1$ and $r=0.8$ for $\Gamma_2$, respectively.
The interface jump conditions $p({\bf x})$, $q({\bf x})$, and the functions $f({\bf x})$, $g({\bf x})$ are set to give
the exact solution of the above model problem~\eqref{eq_poi_jump} as
\begin{equation}\label{ex_jump}
u(x,y) = \left\{
\begin{array}{cc}
-\sin(\pi x) \sin(\pi y) &  \text{in } \Omega_1, \\
e^{-x^2-y^2} &  \text{in } \Omega_2, \\
\sin(\pi x) \sin(\pi y) &  \text{in } \Omega_3.
\end{array}
\right.
\end{equation}
In Figure~\ref{sec33_fig1}, the exact solution plot and an example of sampling points over the subdomain partition
are presented.

According to our model description,
we employ three local neural networks to approximate the model problem solution in each subdomain $\Omega_i$.
The local neural network and sample information are as follows:
\begin{itemize}
\item {\bf Network structure} : 4 hidden layers with 30 neurons per hidden layer in each local neural network.
\item {\bf The number of total parameters (weights, bias) of Network} : 2,911 in each local neural network.
\item {\bf Training samples in $\Omega$} : 2,000 samples for whole interior, 400 samples for four boundary edges (100 per each edge), 150 samples for $\Gamma_1$, and 300 samples for $\Gamma_2$.
\end{itemize}

The error results are reported in Table~\ref{sec33_tb1}.
In Algorithm~\ref{algo3}, the errors become comparable to those
in Algorithms~\ref{algo1} and \ref{algo2} at the later training epochs
though the errors are much larger at the early training epochs.
The augmented Lagrangian term also works effectively for the interface problem
to reduce the errors further in the trained local neural network solutions, see the error results
in Algorithms~\ref{algo1} and \ref{algo2}.

\begin{table}[h!]
\begin{center}
  \caption{
  Error values ${\epsilon}_u$ of the Poisson interface problem in \eqref{eq_poi_jump}-\eqref{ex_jump}
  with three local neural networks:
  The learning rate $\alpha_0 = 0.1$ is used for algorithms $\mathcal{A}_2$ and  $\mathcal{A}_3$.
  $N_l$ is the number of training epochs of local problems in algorithm~$\mathcal{A}_3$ for every outer iteration.
  The numbers are the mean value of errors (and standard deviation) (below: the number of data communication) obtained from five different training sample sets.}
{\footnotesize \renewcommand{\arraystretch}{1.2}
    \begin{tabular}{clllll}
        \Xhline{3\arrayrulewidth}
\multicolumn{2}{c}{\bf Poisson interfaces}  & Max epochs &  &  &  \\
Algorithm & Parameters  & 10,000 epochs & 20,000 epochs & 50,000 epochs & 100,000 epochs
\\[1mm]         \Xhline{0.8\arrayrulewidth}
$\mathcal{A}_1$& &
0.00729 {\scriptsize(1.08e-03)} \vspace{-0.15cm}&
0.00496 {\scriptsize(5.51e-04)}&
0.00328 {\scriptsize(4.60e-04)}&
0.00260 {\scriptsize(4.05e-04)} \\
& & \vspace{-0.01cm} {\scriptsize(10000)}
& \vspace{-0.01cm} {\scriptsize(20000)}
& \vspace{-0.01cm} {\scriptsize(50000)}
& \vspace{-0.01cm} {\scriptsize(100000)}
\\ \Xhline{0.8\arrayrulewidth}
$\mathcal{A}_2$& $\alpha_\lambda=0.1 $&
0.00372  {\scriptsize(8.03e-04)} \vspace{-0.15cm}&
0.00218  {\scriptsize(4.39e-04)}&
0.00128  {\scriptsize(2.41e-04)}&
0.00105  {\scriptsize(8.57e-05)} \\
& & \vspace{-0.01cm} {\scriptsize(10000)}
& \vspace{-0.01cm} {\scriptsize(20000)}
& \vspace{-0.01cm} {\scriptsize(50000)}
& \vspace{-0.01cm} {\scriptsize(100000)}
\\ \Xhline{0.8\arrayrulewidth}
$\mathcal{A}_3$& $N_l = 100, \alpha_\lambda = 0.005$&
0.00823  {\scriptsize(1.04e-03)}  \vspace{-0.15cm}&
0.00467  {\scriptsize(1.28e-03)}&
0.00279  {\scriptsize(5.82e-04)}&
0.00231  {\scriptsize(3.92e-04)} \\
& & \vspace{-0.01cm} {\scriptsize(100)}
& \vspace{-0.01cm} {\scriptsize(200)}
& \vspace{-0.01cm} {\scriptsize(500)}
& \vspace{-0.01cm} {\scriptsize(1000)}
\\
$\mathcal{A}_3$& $N_l = 200, \alpha_\lambda = 0.003$&
0.00846  {\scriptsize(1.21e-03)}  \vspace{-0.15cm}&
0.00439  {\scriptsize(1.98e-03)}&
0.00283  {\scriptsize(3.24e-04)}&
0.00226  {\scriptsize(3.48e-04)} \\
& & \vspace{-0.01cm} {\scriptsize(50)}
& \vspace{-0.01cm} {\scriptsize(100)}
& \vspace{-0.01cm} {\scriptsize(200)}
& \vspace{-0.01cm} {\scriptsize(500)}
\\
$\mathcal{A}_3$& $N_l = 500, \alpha_\lambda = 0.001$&
0.03044  {\scriptsize(1.07e-02)}  \vspace{-0.15cm}&
0.00569  {\scriptsize(1.72e-03)}&
0.00278  {\scriptsize(6.33e-04)}&
0.00220  {\scriptsize(5.19e-04)} \\
& & \vspace{-0.01cm} {\scriptsize(20)}
& \vspace{-0.01cm} {\scriptsize(40)}
& \vspace{-0.01cm} {\scriptsize(100)}
& \vspace{-0.01cm} {\scriptsize(200)}
\\
$\mathcal{A}_3$& $N_l = 1,000, \alpha_\lambda = 0.0005$&
0.14846  {\scriptsize(2.51e-02)}  \vspace{-0.15cm}&
0.04484  {\scriptsize(9.17e-03)}&
0.00337  {\scriptsize(1.06e-03)}&
0.00277  {\scriptsize(1.22e-03)} \\
& & \vspace{-0.01cm} {\scriptsize(10)}
& \vspace{-0.01cm} {\scriptsize(20)}
& \vspace{-0.01cm} {\scriptsize(50)}
& \vspace{-0.01cm} {\scriptsize(100)}
\\
[0.5mm]     \Xhline{3\arrayrulewidth}
    \end{tabular}
}
\label{sec33_tb1}
\vskip-.7truecm
\end{center}
\end{table}

\subsection{Stokes flow with immersed interface}

In this subsection, we consider the following Stokes problem:
\begin{equation}\label{eq_stokes}
\begin{split}
- \mu \triangle \mathbf{u} ( \mathbf{x} )  + \nabla p (\mathbf{x}) & =  \mathbf{h} (\mathbf{x}) + \mathbf{f}( \mathbf{x})  \; \text{ in } \Omega, \\
\nabla \cdot \mathbf{u} (\mathbf{x}) & =  0  \; \text{ in } \Omega, \\
\mathbf{u}({\bf x}) &= \mathbf{g}( {\bf x})  \; \text{ on } \partial \Omega, \\
\end{split}
\end{equation}
where $p(\mathbf{x})$ is the pressure, $\mathbf{u}(\mathbf{x}) := ( u (\mathbf{x}), v(\mathbf{x}) )$ is the fluid velocity in a two dimensional domain $\Omega$, and $\mu$ is a constant viscosity.
The term $\mathbf{f}(\mathbf{x})$ corresponds to the singular force exerted by the immersed interface $\Gamma$ inside the domain $\Omega$ and the term $\mathbf{h}(\mathbf{x})$ to the external force.
Following the immersed boundary method~\cite{peskin2002immersed}, we assume that the singular force $\mathbf{f}(\mathbf{x})$ is supported only along the interface $\Gamma = \{ \mathbf{X} (s) = (X(s), Y(s)),  \ 0 \leq s \leq L \}$ and it can be written as
\begin{equation}\label{f:delta}
 \mathbf{f}(\mathbf{x}) = \int_{\Gamma} \mathbf{F}(\mathbf{X}(s)) \delta (\mathbf{x} - \mathbf{X}(s)) \ d{\mathbf{X}(s)},
\end{equation}
where $\delta$ is the Dirac delta function. The force $\mathbf{F}$ on the interface can be further decomposed into its normal and tangential components,
$$ \mathbf{F} = (\mathbf{F} \cdot \mathbf{n} ) \mathbf{n} + (\mathbf{F} \cdot \mathbf{t})\mathbf{t} = F_{n} \mathbf{n} + F_{\tau} \mathbf{t},$$
where we use the notations $\mathbf{n}$ and $\mathbf{t}$ for the normal and tangential directions on the interface $\Gamma$, respectively. Since the force $\mathbf{f}$ has a delta function singularity along the interface, velocity and the pressure are also nonsmooth across the interface.
By the immersed interface force $\mathbf{F}$, the jump conditions on $\mathbf{u}$ and $p$ can be obtained as the last two equations in \eqref{eq_stokes_if}. Using these jump conditions, the Stokes problem~\eqref{eq_stokes} is rewritten into the following Stokes interface problem on ${\mathbf{u}}$ and $p$:
\begin{equation}\label{eq_stokes_if}
\begin{split}
- \mu \triangle \mathbf{u}  + \nabla p & =  \mathbf{h}  \; \text{ in } \Omega \setminus \Gamma, \\
\nabla \cdot \mathbf{u} & =  0  \; \text{ in } \Omega \setminus \Gamma, \\
\mathbf{u} &= \mathbf{g}  \; \text{ on } \partial \Omega, \\
\left\llbracket  \mathbf{u} \right\rrbracket_{\Gamma} = 0, \quad
\left\llbracket  \mu \frac{\partial \mathbf{u}}{\partial {n}} \right\rrbracket _{\Gamma} &=-F_\tau \mathbf{t} \; \text{ on } \Gamma, \\
\left\llbracket p \right\rrbracket_{\Gamma} = F_{n}, \quad
\left\llbracket  \frac{\partial p}{\partial {n}} \right\rrbracket _{\Gamma} &= \frac{\partial F_\tau}{\partial {\mathbf{X}(s)} \cdot \mathbf{t} +
\left\llbracket \mathbf{h} \right\rrbracket \cdot \mathbf{n} } \; \text{ on } \Gamma,
\end{split}
\end{equation}
where we do not need to deal with the dirac delta function in \eqref{f:delta}
by replacing the immersed interface forcing term $\mathbf{f}$
with the above interface conditions on the solution $\mathbf{u}$ and $p$.

\begin{figure}[htb!]
 \centering
  \includegraphics[height=3.2cm]{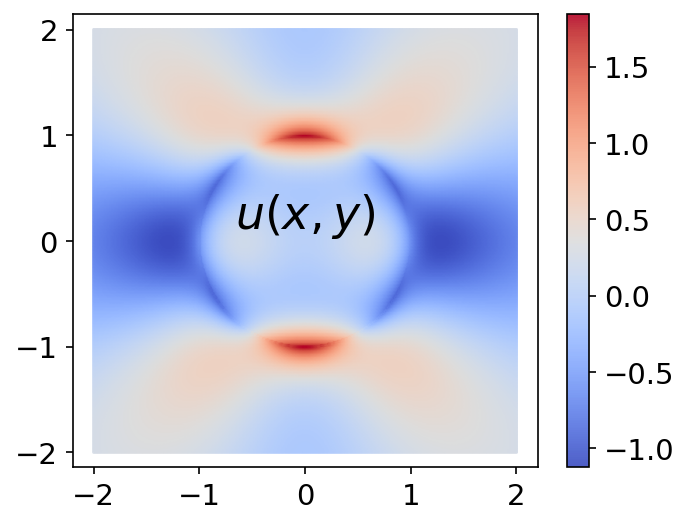}
  \includegraphics[height=3.2cm]{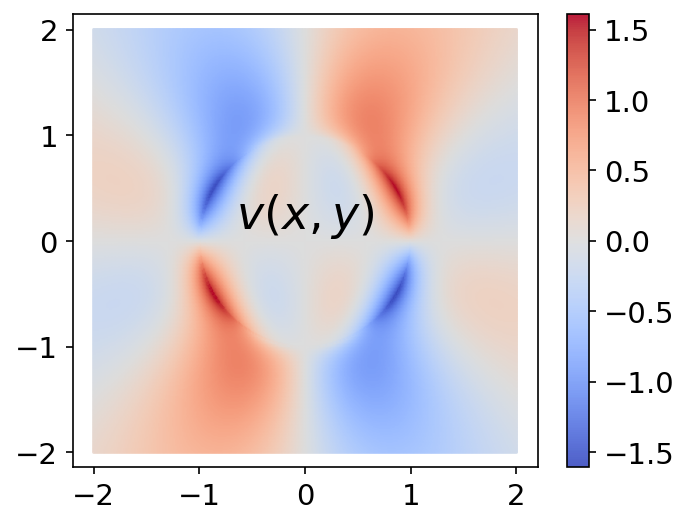}
  \includegraphics[height=3.2cm]{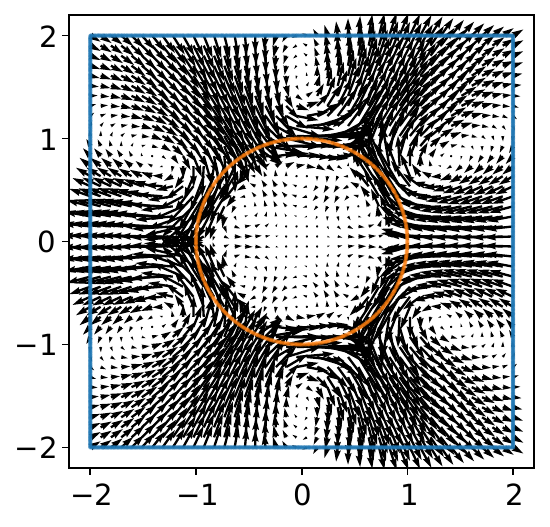}
  \includegraphics[height=3.2cm]{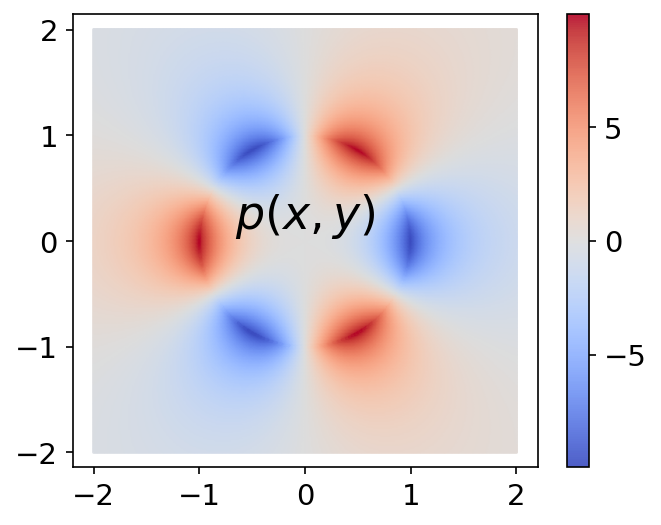}
  \caption{Exact solution plot for \eqref{ex_stokes} and flow quivers: $u(x,y)$ (first), $v(x,y)$ (second), flow quivers (third), $p(x,y)$ (fourth).}
  \label{sec34_fig1}
  \end{figure}

In the following numerical experiment, we set the computational domain $\Omega = (-2,2)^2$, and the interface $\Gamma = \{ \mathbf{X}(s) = ( \cos(s), \sin(s) ) , \ 0 \leq s \leq 2\pi \}$,
where $s$ is the arc length parameter. We also simply set $\mu=1$.
We consider the case in which the immersed force exerts only along the tangential direction and is given by
$$ F_\tau(\mathbf{X}(s)) \mathbf{t}  = 20 \sin (3 s) \mathbf{X}'(s),\quad {F}_n=0.$$
We note that for the above model problem the corresponding exact solution $({\bf u} \,(:=(u,v)),\, p)$ is given as follows:
\begin{equation}\label{ex_stokes}
\begin{split}
u(r,s) &= \left\{
\begin{array}{cc}
\frac{10}{8} r^2 \cos(2s) + \frac{10}{16} r^4 \cos(4s) - \frac{10}{4} r^4 \cos(2 s), &  r<1, \\
-\frac{10}{8} r^{-2} \cos(2s) + \frac{50}{16} r^{-4} \cos(4s) - \frac{10}{4} r^{-2} \cos(4 s), &  r \geq 1,
\end{array}
\right.\\
v(r,s) &= \left\{
\begin{array}{cc}
-\frac{10}{8} r^2 \sin(2s) + \frac{10}{16} r^4 \sin(4s) + \frac{10}{4} r^4 \sin(2 s), &  r<1, \\
\frac{10}{8} r^{-2} \sin(2s) + \frac{50}{16} r^{-4} \sin(4s) - \frac{10}{4} r^{-2} \sin(4 s), &  r \geq 1,
\end{array}
\right.\\
p(r,s) &= \left\{
\begin{array}{cc}
-10r^3 \cos(3s), &  r<1, \\
-10r^{-3} \cos(3s), &  r \geq 1,
\end{array}
\right.
\end{split}
\end{equation}
The exact solution plot and flow quivers are also presented in Figure~\ref{sec34_fig1}.

\begin{figure}[htb!]
 \centering
  \includegraphics[width=17cm]{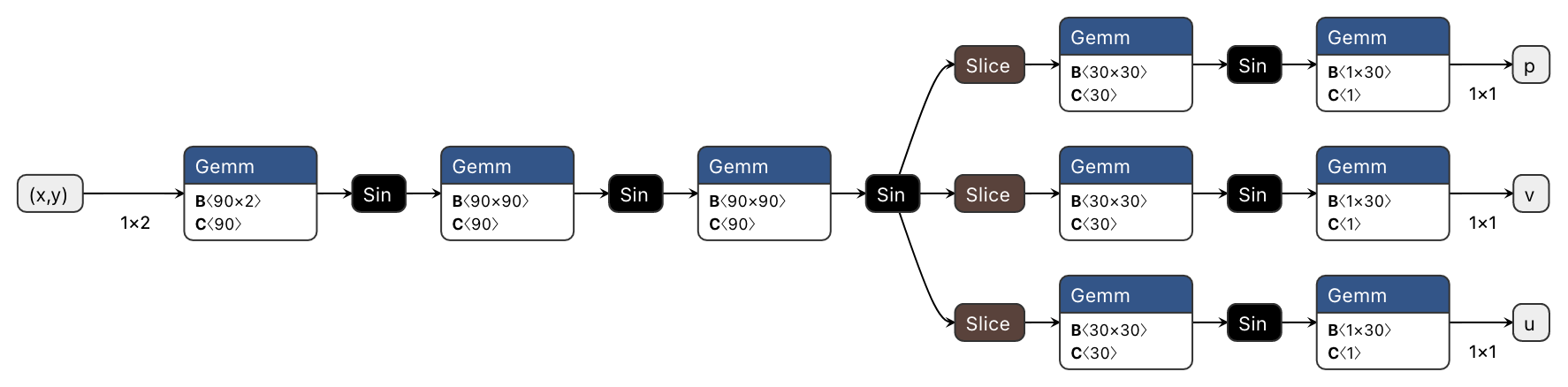}
  \caption{Local network structure for the Stokes problem.}
  \label{sec34_fig2}
  \end{figure}

In our neural network approximation, we employ the following settings for the neural networks and the training samples. We partition the domain into the two subdomains, the interior part of the interface $\Gamma$ and the outside of the interface $\Gamma$. For each subdomain, we set the corresponding local neural network function as described below:
\begin{itemize}
\item {\bf Network structure} : To strengthen the independence of each output variable $u, v$, and $p$, we construct three hidden layers of 90 neurons from the input, and then slice the 90 neurons in the third hidden layer into three branches of 30 neurons. Each sliced branch is followed by a hidden layer of 30 neurons
    and the last output layer. The outputs from the three sliced branches give the values $u, v$, and $p$.
    The resulting network structure is then obtained as seen in Figure~\ref{sec34_fig2}.
\item {\bf The number of total parameters (weights, bias) of Network} : 19,533 for each local neural network.
\item {\bf Training samples} : 5,000 samples for $\Omega$ and 800 samples for boundary (200 for each edge of $\Omega$) and 200 samples for $\Gamma$.
\end{itemize}
For each subdomain $\Omega_i$, we set a single network $\mathcal{N}_i(x; \theta_i)$
as in Figure~\ref{sec34_fig2}, and use them in our computation.
In our numerical experiments, we observed that this kind of network structure
performs better than using independent neural network functions to each unknown, $u$, $v$, and $p$, i.e.,
giving smaller errors in the resulting neural network approximation.
We also note that each single network $\mathcal{N}_i(x;\theta_i)$ has three output values and we denote them by $U_i$, $V_i$ and $P_i$. Below, we use the notation ${\bf U}_i:=(U_i,V_i)$ for the local velocity approximation.

For the Stokes interface problem, we form the following localized loss terms,
\begin{align*}
\mathcal{L}_{f,i}&:=\frac{1}{|X_{\Omega_i}|}\sum_{{\bf x} \in X_{\Omega_i}} | -\mu \triangle {\bf U}_i ({\bf x};\theta_i)+\nabla P_i({\bf x};\theta_i)-{\bf h}({\bf x})|^2,\\
\mathcal{L}_{div,i}&:=\frac{1}{|X_{\Omega_i}|}\sum_{{\bf x} \in X_{\Omega_i}} ( \nabla \cdot {\bf U}_i({\bf x}; \theta_i) )^2,\\
\mathcal{L}_{g,i}&:=\frac{1}{|X_{\Gamma_{i0}}|} \sum_{ {\bf x} \in X_{\Gamma_{i0}}} | {\bf U}_i ({\bf x};\theta_i)-{\bf g}({\bf x})|^2,\\
\mathcal{F}_{u,i}&:=\sum_{j \in s(i)} \frac{1}{|X_{\Gamma_{ij}}|}
 \sum_{{\bf x} \in X_{\Gamma_{ij}}} \left( |{\bf U}_i ({\bf x};\theta_i)-{\bf \tilde{U}}_{ij} ({\bf x}) |^2
 + (P_i({\bf x};\theta_i)- \tilde{P}_{ij}({\bf x}))^2\right),\\
\mathcal{F}_{n,i}&:=\sum_{j \in s(i)} \frac{1}{|X_{\Gamma_{ij}}|}
 \sum_{{\bf x} \in X_{\Gamma_{ij}}} \left( \left|\mu \frac{\partial {\bf U}_i}{\partial n} ({\bf x};\theta_i)-{\bf \tilde{U}}_{n,ij} ({\bf x}) \right|^2
 + \left(\frac{\partial P_i}{\partial n} ({\bf x};\theta_i)- \tilde{P}_{n,ij}({\bf x}) \right)^2\right),
\end{align*}
where ${\bf \tilde{U}}_{ij} ({\bf x})$, $\tilde{P}_{ij}({\bf x})$, ${\bf \tilde{U}}_{n,ij} ({\bf x})$, and $\tilde{P}_{n,ij}({\bf x})$ are defined similarly as in \eqref{u_tildeij} and \eqref{u_tildeji}.
For the Stokes interface problem, we introduce three kinds of Lagrange multipliers;
First, ${\lambda}_{u,i0}$ and $\lambda_{v,i0}$ are for the boundary condition of velocity variable ${\bf u}$, i.e.,
$$\Lambda_{i0}(\theta_i,(\lambda_{u,i0},\lambda_{v,i0})):=\sum_{{\bf x} \in X(\Gamma_{i0})} (\lambda_{u,i0}({\bf x}),\lambda_{v,i0}({\bf x}))\cdot
( {\bf U}_i ({\bf x};\theta_i)-{\bf g} ({\bf x}) ).$$
Second, $\lambda_{u,ij}$, $\lambda_{v,ij}$ and $\lambda_{p,ij}$ are for the interface conditions on
the velocity and pressure,
\begin{equation*}
\begin{split}
&\Lambda_{i}(\theta_i,(\lambda_{u,ij},\lambda_{v,ij},\lambda_{p,ij})_j)\\
&:=\sum_{j \in s(i)} \sum_{{\bf x} \in X(\Gamma_{ij})} \left( (\lambda_{u,ij}({\bf x}),\lambda_{v,ij}({\bf x}))\cdot
( {\bf U}_i ({\bf x};\theta_i)-{\bf \tilde{U}}_{ij} ({\bf x}))+ \lambda_{p,ij}({\bf x}) ( P_i({\bf x};\theta_i)- \tilde{P}_{ij}({\bf x}))\right).
\end{split}
\end{equation*}
In addition, we introduce vector unknowns $\lambda_{div,i}$ such that their length are identical to the number of training sample points in $\Omega_i$.
These Lagrange multipliers are introduced
for the divergence-free constraint on the velocity approximation, $\nabla \cdot \mathbf{u} = 0$,
\begin{equation}\label{lag_div}
\Lambda_{{div},i} =  \sum_{{\bf x} \in X_{ \Omega_i }} \lambda_{div,i}({\bf x}) \nabla \cdot {\bf U }_i ({\bf x};\theta_i),
\end{equation}
and $\lambda_{div,i}$ are updated by the gradient ascent method,
\begin{equation}
\lambda_{div,i}^{(\ell)} ({\bf x}) = \lambda_{div,i}^{(\ell-1)} ({\bf x})+ \alpha_d  \nabla \cdot {\bf U}_i ({\bf x};\theta_i^{(\ell)})
\end{equation}
with a suitable learning rate $\alpha_d$.
We also observed that the augmented Lagrangian term for the divergence free condition further increases
the training performance in our neural network approximation by enforcing the physical property of
the solution much stronger.
Lagrange multipliers $\lambda_{u,i0}$ and $\lambda_{v,i0}$, for the boundary condition, are initialized to zero vectors once before the training start.
On the other hand, Lagrange multipliers, $\lambda_{u,ij}$, $\lambda_{v,ij}$, $\lambda_{p,ij}$ and $\lambda_{div,i}$ are initialized as the gradient of the loss function before staring a new outer iteration in Algorithm~\ref{algo3}.
With all the previously defined individual loss terms and the augmented Lagrangian terms,
for the Stokes flow with the immerse interface, in Algorithm~\ref{algo1} we employed the following localized loss function,
$$\mathcal{J}_i(\theta_i):=\mathcal{L}_{f,i}+\mathcal{L}_{div,i}+\mathcal{L}_{g,i}+\mathcal{F}_{u,i}+\mathcal{F}_{n,i},$$
and in Algorithms~\ref{algo2} and \ref{algo3}  we employed the localized loss function enriched with the augmented Lagrangian terms,
\begin{equation*}
\begin{split}
&\mathcal{J}_{i,\Lambda}(\theta_i,(\lambda_{u,i0},\lambda_{v,i0}), (\lambda_{u,ij},\lambda_{v,ij}, \lambda_{p,ij})_j, \lambda_{div,i})\\
&:=\mathcal{J}_i(\theta_i)+\Lambda_{i0}(\theta_i,(\lambda_{u,i0},\lambda_{v,i0}))+\Lambda_i(\theta_i,(\lambda_{u,ij},\lambda_{v,ij}, \lambda_{p,ij})_j)+ \Lambda_{div,i}(\theta_i,\lambda_{div,i}).
\end{split}
\end{equation*}

For the trained neural network approximate solution, we compute the relative $L^2$-errors for the velocity and pressure as
\begin{equation}\label{Stokes_err}
{\epsilon}_{\mathbf{u}}(X) = \sqrt{\frac{ \sum_{{\bf x} \in X} | {\bf u}({\bf x}) - {\bf U} ({\bf x}; \theta^*) |^2}{ \sum_{{\bf x}\in X} |{\bf u} ({\bf x})|^2 } }, \quad {\epsilon}_{p}(X) = \sqrt{ \frac{ \sum_{{\bf x} \in X} \left(  p({\bf x}) - P ({\bf x}; \theta^*) \right)^2 }{ \sum_{{\bf x}\in X} p({\bf x})^2 } },
\end{equation}
where
$${\bf U}({\bf x};\theta^*)={\bf U}_i({\bf x};\theta_i^*),\; P({\bf x};\theta^*)=P_i({\bf x};\theta_i^*),\;\text{ for } \; {\bf x} \in \Omega_i.$$
The error results are reported in Table~\ref{sec34_tb1}. For the Stokes immersed-interface problem,
Algorithm~\ref{algo2} outperforms Algorithm~\ref{algo1}, showing that the augmented Lagrangian terms
help to reduce errors in the neural network approximation.
Algorithm~\ref{algo3}, the iteration version of Algorithm~\ref{algo2}, gives convergent approximate solutions
with acceptable error results as proceeding more training epochs.
To enhance the training performance further in the iteration method,
acceleration algorithms for the iteration method~$\mathcal{A}_3$ need to be investigated by
using well-developed preconditioning tools in the classical domain decomposition methods~\cite{TW-Book}.

\begin{table}[h!]
\begin{center}
  \caption{
   Error values ${\epsilon}_{\mathbf{u}}$ and ${\epsilon}_{p}$ of the Stokes problem~\eqref{eq_stokes_if} with the immersed-interface singular force example~\eqref{ex_stokes}. The domain is partitioned into two subdomains and two local neural networks are employed:
   The learning rate $\alpha_0 = 0.1$ is used for algorithms $\mathcal{A}_2$ and  $\mathcal{A}_3$.
  $N_l$ is the number of training epochs of local problems in algorithm~$\mathcal{A}_3$ for every outer iteration.
  The numbers are the mean value of errors (and standard deviation) (below: the number of data communication) obtained from five different training sample sets.}
\label{sec34_tb1}
{\footnotesize \renewcommand{\arraystretch}{1.2}
    \begin{tabular}{clllll}
        \Xhline{3\arrayrulewidth}
\multicolumn{2}{c}{\bf Stokes interface}   & Max epochs &  &  &  \\
Algorithm & Parameters  & 10,000 epochs & 20,000 epochs & 50,000 epochs & 100,000 epochs \\
& & ${\epsilon}_{\mathbf{u}}$ \hspace{0.65cm} ${\epsilon}_p$
& ${\epsilon}_{\mathbf{u}}$ \hspace{0.65cm} ${\epsilon}_p$
& ${\epsilon}_{\mathbf{u}}$ \hspace{0.65cm} ${\epsilon}_p$
& ${\epsilon}_{\mathbf{u}}$ \hspace{0.65cm} ${\epsilon}_p$
\\[1mm]         \Xhline{0.8\arrayrulewidth}
$\mathcal{A}_1$& &
0.0206, \ 0.0317 \vspace{-0.15cm}&
0.0111, \ 0.0190&
0.0057, \ 0.0105&
0.0036, \ 0.0070 \\
& & \vspace{-0.01cm} {\scriptsize(10000)}
& \vspace{-0.01cm} {\scriptsize(20000)}
& \vspace{-0.01cm} {\scriptsize(50000)}
& \vspace{-0.01cm} {\scriptsize(100000)}
\\ \Xhline{0.8\arrayrulewidth}
$\mathcal{A}_2$& $\alpha_\lambda = 0.1$, $\alpha_d = 0.1$&
0.0056, \ 0.0105 \vspace{-0.15cm}&
0.0028, \ 0.0042&
0.0018, \ 0.0020&
0.0012, \ 0.0013 \\
& & \vspace{-0.01cm} {\scriptsize(10000)}
& \vspace{-0.01cm} {\scriptsize(20000)}
& \vspace{-0.01cm} {\scriptsize(50000)}
& \vspace{-0.01cm} {\scriptsize(100000)}
\\ \Xhline{0.8\arrayrulewidth}
$\mathcal{A}_3$& $N_l = 100, \alpha_\lambda = 0.02$, $\alpha_d=0.2$&
0.0136, \ 0.0199 \vspace{-0.15cm}&
0.0115, \ 0.0147&
0.0081, \ 0.0107&
0.0068, \ 0.0076 \\
& & \vspace{-0.01cm} {\scriptsize(100)}
& \vspace{-0.01cm} {\scriptsize(200)}
& \vspace{-0.01cm} {\scriptsize(500)}
& \vspace{-0.01cm} {\scriptsize(1000)}
\\
$\mathcal{A}_3$& $N_l = 200, \alpha_\lambda = 0.001$, $\alpha_d=0.1$&
0.0132, \ 0.0151 \vspace{-0.15cm}&
0.0073, \ 0.0100&
0.0066, \ 0.0080&
0.0070, \ 0.0064 \\
& & \vspace{-0.01cm} {\scriptsize(50)}
& \vspace{-0.01cm} {\scriptsize(100)}
& \vspace{-0.01cm} {\scriptsize(200)}
& \vspace{-0.01cm} {\scriptsize(500)}
\\
$\mathcal{A}_3$& $N_l = 500, \alpha_\lambda = 0$, $\alpha_d=0.05$&
0.1700, \ 0.1246 \vspace{-0.15cm}&
0.0375, \ 0.0288&
0.0058, \ 0.0058&
0.0049, \ 0.0051 \\
& & \vspace{-0.01cm} {\scriptsize(20)}
& \vspace{-0.01cm} {\scriptsize(40)}
& \vspace{-0.01cm} {\scriptsize(100)}
& \vspace{-0.01cm} {\scriptsize(200)}
\\
$\mathcal{A}_3$& $N_l = 1,000, \alpha_\lambda = 0$, $\alpha_d=0.01$&
0.3843, \ 0.3022 \vspace{-0.15cm}&
0.2121, \ 0.1630&
0.0343, \ 0.0272&
0.0071, \ 0.0078 \\
& & \vspace{-0.01cm} {\scriptsize(10)}
& \vspace{-0.01cm} {\scriptsize(20)}
& \vspace{-0.01cm} {\scriptsize(50)}
& \vspace{-0.01cm} {\scriptsize(100)}
\\
[0.5mm]     \Xhline{3\arrayrulewidth}
    \end{tabular}
}
\vskip-.7truecm
\end{center}
\end{table}

\section{Conclusions}\label{sec:conclude}
Partitioned neural network approximation to PDEs is proposed for non-overlapping subdomain partitions of the problem domain. Each local neural network approximates the solution in each subdomain. Localized loss functions and augmented Lagrangian terms are introduced to enhance the neural network solution accuracy and the training performance. To reduce communication cost among the neighboring local neural networks during the parameter training epochs, an iteration algorithm, Algorithm 3, is also proposed and its convergence and training performance are numerically studied. Numerical results also confirm the capability of the proposed methods for application problems with discontinuous coefficients, interface discontinuity, and an immersed interface.

For the iteration algorithm, Algorithm 3, its convergence needs to be theoretically analyzed and improved to enhance the parallel computing scalability. Algorithm 3 is similar to the FETI algorithms~\cite{FETI,local-FETI,total-FETI}
and our outlook is that its convergence analysis can be done by borrowing various analytical tools in the FETI algorithms.
Our future research will be focused on the theoretical convergence analysis of Algorithm 3
and its convergence enhancement by utilizing many successful preconditioning techniques
in non-overlapping domain decomposition algorithms~\cite{TW-Book}.

\bibliographystyle{elsarticle-num}
\bibliography{ddforpinn}

\begin{thebibliography}{10}
\expandafter\ifx\csname url\endcsname\relax
  \def\url#1{\texttt{#1}}\fi
\expandafter\ifx\csname urlprefix\endcsname\relax\def\urlprefix{URL }\fi
\expandafter\ifx\csname href\endcsname\relax
  \def\href#1#2{#2} \def\path#1{#1}\fi

\bibitem{sirignano2018}
J.~Sirignano, K.~Spiliopoulos, D{GM}: a deep learning algorithm for solving
  partial differential equations, J. Comput. Phys. 375 (2018) 1339--1364.

\bibitem{yu2018deep}
B.~Yu, et~al., The deep {R}itz method: a deep learning-based numerical
  algorithm for solving variational problems, Communications in Mathematics and
  Statistics 6~(1) (2018) 1--12.

\bibitem{raissi2019}
M.~Raissi, P.~Perdikaris, G.~E. Karniadakis, Physics-informed neural networks:
  a deep learning framework for solving forward and inverse problems involving
  nonlinear partial differential equations, J. Comput. Phys. 378 (2019)
  686--707.

\bibitem{cPINN}
A.~D. Jagtap, E.~Kharazmi, G.~E. Karniadakis, Conservative physics-informed
  neural networks on discrete domains for conservation laws: {A}pplications to
  forward and inverse problems, Computer Methods in Applied Mechanics and
  Engineering 365 (2020) 113028.

\bibitem{xPINN}
A.~D. Jagtap, G.~E. Karniadakis, Extended physics-informed neural networks
  (xpinns): A generalized space-time domain decomposition based deep learning
  framework for nonlinear partial differential equations, Communications in
  Computational Physics 28~(5) (2020) 2002--2041.

\bibitem{FBPINN}
B.~Moseley, A.~Markham, T.~Nissen-Meyer, {F}inite {B}asis {P}hysics-{I}nformed
  {N}eural {N}etworks ({FBPINN}s): a scalable domain decomposition approach for
  solving differential equations, arXiv preprint arXiv:2107.07871 (2021).

\bibitem{TW-Book}
A.~Toselli, O.~Widlund, Domain decomposition methods---algorithms and theory,
  Vol.~34 of Springer Series in Computational Mathematics, Springer-Verlag,
  Berlin, 2005.

\bibitem{DD26-kim-yang}
H.~H. Kim, H.~J. Yang, Domain decomposition algorithms for physics-informed
  neural networks, in: Proceedings of the 26th International Conference on
  Domain Decomposition Methods, 2021.

\bibitem{DD27-kim-yang}
H.~H. Kim, H.~J. Yang, Domain decomposition algorithms for neural network
  approximation of partial differential equations, in: Proceedings of the 27th
  International Conference on Domain Decomposition Methods (To appear).

\bibitem{arXiv-yang-kim-2023}
H.~J. Yang, H.~H. Kim, Iterative algorithms for partitioned neural network
  approximation to partial differential equations, Submitted (2023).

\bibitem{li2019}
K.~Li, K.~Tang, T.~Wu, Q.~Liao, {D3M}: {A} deep domain decomposition method for
  partial differential equations, IEEE Access 8 (2019) 5283--5294.

\bibitem{li2020-pro}
W.~Li, X.~Xiang, Y.~Xu, Deep domain decomposition method: {E}lliptic problems,
  in: Mathematical and Scientific Machine Learning, PMLR, 2020, pp. 269--286.

\bibitem{Son2023}
H.~Son, S.~W. Cho, H.~J. Hwang, Enhanced physics-informed neural networks with
  augmented {L}agrangian relaxation method ({AL-PINN}s), Neurocomputing (2023)
  126424.

\bibitem{Huang2021}
J.~Huang, H.~Wang, T.~Zhou, An augmented {L}agrangian deep learning method for
  variational problems with essential boundary conditions, arXiv preprint
  arXiv:2106.14348 (2021).

\bibitem{Scalable-PINN2022}
K.~Shukla, M.~Xu, N.~Trask, G.~E. Karniadakis, Scalable algorithms for
  physics-informed neural and graph networks, Data-Centric Engineering 3 (2022)
  e24.

\bibitem{Substruct}
A.~Klawonn, O.~Widlund, {FETI} and {N}eumann-{N}eumann iterative substructuring
  methods: connections and new results, Communications on Pure and Applied
  Mathematics: A Journal Issued by the Courant Institute of Mathematical
  Sciences 54~(1) (2001) 57--90.

\bibitem{FETI}
C.~Farhat, F.-X. Roux, A method of finite element tearing and interconnecting
  and its parallel solution algorithm, International journal for numerical
  methods in engineering 32~(6) (1991) 1205--1227.

\bibitem{BDD}
J.~Mandel, Balancing domain decomposition, Communications in numerical methods
  in engineering 9~(3) (1993) 233--241.

\bibitem{local-FETI}
K.~Park, C.~A. Felippa, U.~Gumaste, A localized version of the method of
  {L}agrange multipliers and its applications, Computational Mechanics 24~(6)
  (2000) 476--490.

\bibitem{total-FETI}
Z.~Dost{\'a}l, D.~Hor{\'a}k, R.~Ku{\v{c}}era, Total {FETI}:an easier
  implementable variant of the {FETI} method for numerical solution of elliptic
  {PDE}, Communications in Numerical Methods in Engineering 22~(12) (2006)
  1155--1162.

\bibitem{sitzmann2020implicit}
V.~Sitzmann, J.~Martel, A.~Bergman, D.~Lindell, G.~Wetzstein, Implicit neural
  representations with periodic activation functions, Advances in neural
  information processing systems 33 (2020) 7462--7473.

\bibitem{adam}
D.~P. Kingma, J.~Ba, Adam: A method for stochastic optimization, in:
  International Conference on Learning Representations (ICLR), 2015.

\bibitem{peskin2002immersed}
C.~S. Peskin, The immersed boundary method, Acta numerica 11 (2002) 479--517.

\end{thebibliography}

\end{document}